\documentclass[12pt,reqno]{amsart}

\usepackage{setspace,  
hyperref, epigraph} 

\usepackage{breakurl}   

\usepackage{enumitem} 

\numberwithin{equation}{section}

\numberwithin{figure}{section} 

\DeclareMathOperator{\R}{{\mathbb R}}

\DeclareMathOperator{\N}{{\mathbb N}}

\newcommand{\st}{\textbf{st}}

\newcommand\law{SIA}

\newcommand\Xlone{{}_1\hskip-1.1pt{X}}

\author[J. Bair]{Jacques Bair}\address{J. Bair, HEC-ULG, University of 
Liege, 4000 Belgium}\email{J.Bair@ULiege.be}

\author[P. B{\l}aszczyk]{Piotr B\l{}aszczyk}\address{P. B{\l}aszczyk,
Institute of Mathematics, Pedagogical University of Cracow,
Poland}\email{pb@up.krakow.pl}

\author[R. Ely]{Robert Ely}\address{R. Ely, Department of Mathematics,
University of Idaho, Moscow, ID 83844 US}\email{ely@uidaho.edu}

\author[M. Katz]{Mikhail G. Katz}\address{M. Katz, Department of
Mathematics, Bar Ilan University, Ramat Gan 52900
Israel}\email{katzmik@macs.biu.ac.il}

\author[K. Kuhlemann]{Karl Kuhlemann}\address{K. Kuhlemann, Gottfried
Wilhelm Leibniz University Hannover, D-30167 Hannover,
Germany}\email{kus.kuhlemann@t-online.de}


\begin{document}


\thispagestyle{empty}


\title [Procedures of Leibnizian infinitesimal calculus] 
{Procedures of Leibnizian infinitesimal calculus: An account in three
modern frameworks}

\begin{abstract}
Recent Leibniz scholarship has sought to gauge which foundational
framework provides the most successful account of the procedures of
the Leibnizian calculus (LC).  While many scholars (e.g., Ishiguro,
Levey) opt for a default Weierstrassian framework, Arthur compares LC
to a non-Archimedean framework SIA (Smooth Infinitesimal Analysis) of
Lawvere--Kock--Bell.  We analyze Arthur's comparison and find it rife
with equivocations and misunderstandings on issues including the
non-punctiform nature of the continuum, infi\-nite-sided polygons, and
the fictionality of infinitesimals.  Rabouin and Arthur claim that
Leibniz considers infinities as contradictory, and that Leibniz'
definition of incomparables should be understood as nominal rather
than as semantic.  However, such claims hinge upon a conflation of
Leibnizian notions of bounded infinity and unbounded infinity, a
distinction emphasized by early Knobloch.


The most faithful account of LC is arguably provided by Robinson's
framework.  We exploit an axiomatic framework for infinitesimal
analysis called SPOT (conservative over ZF) to provide a formalisation
of LC, including the bounded/unbounded dichotomy, the
assignable/inassignable dichotomy, the generalized relation of
equality up to negligible terms, and the law of continuity.

\end{abstract}


\maketitle
\tableofcontents

\epigraph{We find here\, [a]\, very important feature of Leibniz's
strategy when dealing with infinite and infinitely small quantities:
nowhere does he express his conviction that they are contradictory
notions.  ---\;Rabouin {\&} Arthur, 2020}

\epigraph{The most popular grounds for denying the existence of
infinitesimals, besides a lingering sense of their simply being
conceptually ``repugnant,'' are (1) the lack of any useful
mathematical application, and\;(2)\;their incompatibility with the
so-called axiom of Archimedes.  --\;Levey, 1998}

\epigraph{[Leibniz] has defined {\ldots} what it is for one quantity
(say,~$dx$) to be incomparable \emph{in relation to} another
(say,~$x$).  Formally, Leibniz's definition gives:

For a given~$x$ and~$dx$,\,~$dx$ INC~$x$ iff~$\neg(\exists n)
(ndx>x)$.  ---\;Rabouin {\&} Arthur, 2020}

\section{Introduction}
\label{one}

Recent Leibniz scholarship has sought to gauge which foundational
framework furnishes the most successful account of the procedures of
the Leibnizian calculus (LC).  While many scholars (e.g., Ishiguro,
Levey) opt for a default Weierstrassian framework, some recent work
explores additional possibilities, such as the non-Archimedean systems
provided by Lawvere's framework and by Robinson's framework.  Siegmund
Probst observes:
\begin{quote}
For us, the question is whether Leibniz based his metaphysical
foundations of the calculus wholly on the concept of the
syncategorematic infinite or whether he pursued also an alternative
approach of accepting infinitesimals similar to what is done in modern
non-standard analysis.  (Probst \cite{Pr18}, 2018, p.\;221)
\end{quote}
We will extend the scope of Probst's question to include the
additional possibility of Lawvere's framework.  We will give summaries
of all three accounts in this introduction.

\subsection
{Grafting of the \emph{Epsilontik} on the calculus of Leibniz}
\label{s11}

The title of this section is inspired by that of Cajori (\cite{Ca23},
1923).  Commenting on the role of the theory of limits, Cajori noted:
``[Leibniz] never actually founded his calculus upon that theory''
\cite[p.\;223]{Ca23}.  Cajori saw this aspect of LC as a shortcoming:
\begin{quote}
At that time both [Leibniz and Newton] used infinitely small
quantities which were dropped when comparatively small.%
\footnote{See note~\ref{f21b} for Boyer's summary of the procedure of
dropping negligible terms, and Section~\ref{s35} for a formalisation.}
It is one of the curiosities in the history of mathematics that this
rough procedure was adopted, even though before this time theories of
\emph{limits} had been worked out in geometry by Giovanni Benedetti,
S.  Stevin, L.\;Valerius, Gregory St.\;Vincent and Tacquet, and in
higher arithmetic by John Wallis.  (ibid.; emphasis added).
\end{quote}
Cajori went on to summarize his view concerning the role of limits as
follows: ``[W]e see on the European continent the grafting of the
Newtonian and pre-Newtonian concepts of limits upon the calculus of
Leibniz'' \cite[p.\;234]{Ca23}.  Cajori was writing a number of years
before the following mathematical developments took place:
\begin{enumerate}
\item
Skolem developed proper extensions of the naturals satisfying the
axioms of Peano Arithmetic in (\cite{Sk33}, 1933 and \cite{Sk34},
1934);
\item
Hewitt introduced hyper-real ideals in (\cite{He48}, 1948);
\item
\L o\'s proved his theorem in (\cite{Lo55}, 1955);
\item
Robinson pioneered his framework for infinitesimal analysis in
(\cite{Ro61}, 1961);
\item
Lawvere pioneered his category-theoretic framework in the 1960s and
70s (see e.g., \cite{La80}).
\end{enumerate}
These developments are not to be construed as breakthroughs in Leibniz
scholarship, but rather as mathematical developments that enabled
modern theories of infinitesimals, one of which was compared to LC by
Arthur (see Section~\ref{s17}).

\subsection{From Berkeley's ghosts to Guicciardini's limits}

In 1734, George Berkeley claimed that infinitesimals were logically
contradictory.  Namely, he claimed that they are simultaneously zero
and nonzero.  There was a bit of a flurry in England at the time
around Berkeley's pamphlet; see e.g., Andersen (\cite{An11}, 2011).
However, it was mostly ignored until the end of the 19th century.  In
the closing decades of the 19th century, the trend-setters in the
mathematics community ruled against infinitesimals, and Berkeley's
criticism became popular again.  This development was due to the fact
that some working mathematicians (though by no means all) agreed with
its \emph{conclusion} (``infinitesimals are contradictory"), not
because they actually analyzed his argument.  Since then, many
historians have quoted Berkeley's claims without analyzing them
properly.  

Robinson's work (\cite{Ro66}, 1966) cast doubt upon Berkeley's
conclusion, leading scholars like David Sherry (\cite{Sh87}, 1987) to
take a fresh look at Berkeley's criticism to see if it actually made
sense, and they were not impressed.  Contrary to Berkeley's claim,
Leibniz never asserted that an infinitesimal literally equals zero.%
\footnote{Possibly Johann Bernoulli did, but there may have been
significant differences between Leibniz and Bernoulli; see Nagel
(\cite{Na08}, 2008).}
On the contrary, Leibniz emphasized repeatedly that he worked with a
generalized relation of equality ``up to" a negligible term (see
Section~\ref{s28}).  In this way, Leibniz justified moves such as the
passage from~$\frac{2x+dx}{a}$ to~$\frac{2x}{a}$ in the calculation
of~$\frac{dy}{dx}$ when~$ay=x^2$ (see Section~\ref{s42}).  Even such a
generalized relation had antecedents.  Thus, Fermat's method of
adequality involved \emph{discarding} the remaining terms in~$E$
(rather than setting them equal to zero); see Str{\o}mholm
(\cite{St68}, 1968, p.\;51).%
\footnote{\label{f21b}As observed by Boyer, ``The general principle
that in an equation involving infinitesimals those of higher order are
to be discarded inasmuch they have no effect on the final result is
sometimes regarded as the basic principle of the differential
calculus. {\ldots} this type of doctrine constituted the central theme
in the developments leading to the calculus of Newton and Leibniz
{\ldots} It was, indeed, precisely upon this general premise that
Leibniz sought to establish the method of differentials''
\cite[p.\;79]{Bo41}.  For a possible formalisation of the procedure of
discarding terms see Section~\ref{s35}.}

Somewhat paradoxically, Robinson's insights appear to have led to an
intensification of the efforts by some historians to demonstrate that
(contrary to Cajori's claim) Leibniz did found his calculus on limits.
Thus, Guicciardini writes:
\begin{quote}
In several of his mature writings, [Leibniz] stated that his
differentials are `well founded,' since they are symbolic
abbreviations for \emph{limit procedures}.  (Guicciardini \cite{Gu18},
2018, p.\;90; emphasis added)
\end{quote}
In a similar vein, Guicciardini claims that 
\begin{quote}
Leibniz distinguished between the search for rigor and the
effectiveness of algorithms. The former was the purview of
metaphysicians, because it ultimately rests on metaphysical
principles, most notably on the `principles of continuity,' on which
Leibniz `well founded,' as he used to say, his method of
\emph{infinitesimals as limits}.  (ibid.; emphasis added)
\end{quote}
Here Guicciardini alludes to Leibnizian references to infinitesimals
as \emph{well-founded fictions}.  However, Leibniz never wrote
anything about any ``method of infinitesimals as limits.''%
\footnote{The view of \emph{infinitesimals as limits} is no late
stratagem for Guicciardini.  Nearly two decades earlier, he already
claimed that ``Leibniz carefully defines the infinitesimal and the
infinite in terms of limit procedures'' (Guicciardini \cite{Gu00},
2000).  The claim occurs in a review of Knobloch (\cite{Kn99}, 1999).
Note that the article under review does not mention the word ``limit''
in connection with the Leibnizian calculus.}
If, as per Guicciardini, Leibniz viewed infinitesimals as limits, then
Leibniz's vision of the infinitesimal calculus was significantly
different from l'Hospital's.  Unlike Guicciardini, editors Bradley et
al.\;of an English edition of l'Hospital's \emph{Analyse des
infiniment petits} take for granted a \emph{continuity} between
Leibniz's and l'Hospital's visions of the calculus:
\begin{quote}
Leibniz' differential calculus tells us how to find the relations
among infinitely small increments~$dx$,~$dy$, etc., among the
variables~$x$,~$y$, etc., in an equation.  L'H\^opital gives the
definition on p.\;2: ``The infinitely small portion by which a
variable quantity continually increases or decreases is called the
\emph{Differential}.''  (Bradley et al.\;\cite{Br15}, 2015, p.\;xvii;
emphasis in the original)
\end{quote}
While there were differences between the views of Leibniz and
l'Hospital (notably with regard to the question of the reality of
infinitesimals), the assumption of a \emph{dis}continuity between them
is a thesis that needs to be argued rather than postulated%
\footnote{\label{f5b}Similar remarks apply to the question of
continuity between Leibniz and other Leibnizians such as Jacob
Hermann; see note~\ref{f16}.  Jesseph analyzes the use of
infinitesimals by Torricelli and Roberval, and concludes: ``By taking
indivisibles as infinitely small magnitudes of the same dimension as
the lines, figures, or solids they compose, [Roberval] could avoid the
paradoxes that seemed to threaten Cavalieri's methods, just as
Torricelli had done'' (\cite{Je21}, 2021, p.\;117).  Thus Torricelli
and Roberval viewed infinitesimals as homogeneous with ordinary
quantities but incomparable with them.  In Section~\ref{V4} we argue
that this was Leibniz's view, as well.  Guicciardini's limit-centered
interpretation of Leibnizian infinitesimals anachronistically
postulates a discontinuity between Leibniz with both the preceding and
the following generations of mathematicians, and is symptomatic of
butterfly-model thinking (see Section~\ref{s17b}).}
(for a discussion of the issue of historical continuity see Katz
\cite{20z}, 2020).  The assumption that Leibniz based his calculus on
limits is akin to Ishiguro's interpretation; see Section~\ref{s13}.

We will examine several modern interpretations of LC, with special
attention to Arthur's comparison of LC to Smooth Infinitesimal
Analysis (SIA); see Section~\ref{s17}.  An introduction to SIA
appears in Section~\ref{s3}.

\subsection{Modern interpretations}
\label{s13}

Interpretations of LC by Bos \cite{Bo74} and Mancosu \cite{Ma89}
contrast with interpretations by Ishiguro (\cite{Is90}, 1990,
Chapter\;5) and Arthur.  These interpretations were explored in
Katz--Sherry (\cite{12e}, 2012 and \cite{13f}, 2013), Bascelli et
al.\;(\cite{16a}, 2016), B{\l}aszczyk et al.\;(\cite{17d}, 2017), Bair
et al.\;(\cite{17b}, 2017 and \cite{18a}, 2018).

Ishiguro interprets Leibnizian infinitesimals as follows:
\begin{quote}
It seems that when we make reference to infinitesimals in a
proposition, we are not designating a fixed magnitude incomparably
smaller than our ordinary magnitudes.  Leibniz is saying that whatever
small magnitude an opponent may present, one can assert the existence
of a smaller magnitude. In other words, we can paraphrase the
proposition with a \emph{universal} proposition with an embedded
\emph{existential} claim.  (Ishiguro \cite{Is90}, 1990, p.\;87;
emphasis added)
\end{quote}
Ishiguro posits that when Leibniz wrote that his inassignable~$dx$, or
alternatively~$\epsilon$, was smaller than every given quantity~$Q$,
what he really meant was a quantifier statement to the effect that for
each given~$Q>0$ there exists an~$\epsilon>0$ such
that~\mbox{$\epsilon<Q$} (or that the error is smaller than~$Q$ if the
increment is smaller than~$\epsilon$).%
\footnote{Thus Ishiguro's Chapter\;5 fits in a modern tradition of
scholarship in search of, to paraphrase Berkeley, the ghosts of
departed \emph{quantifiers} that targets not only Leibniz but also
Cauchy; see e.g., Bair et al.\;(\cite{17a}, 2017).}
This is the standard \emph{syncategorematic} interpretation, via
alternating logical quantifiers, of terms involving infinity.
According to such interpretations, Leibnizian infinitesimals are
\emph{logical fictions}.  Sherry--Katz \cite{14c} contrasted such
interpretations with a \emph{pure fictionalist} interpretation.  Every
occurrence of a logical fiction is eliminable by a quantifier
paraphrase, whereas pure fictions are not.

Arthur has endorsed Ishiguro's reading in a series of articles and
books; see e.g., \cite{Ar07}--\cite{Ar18}, \cite{Ra20}.  The so-called
``Leibniz's syncategorematic infinitesimals'' currently come in two
installments: Arthur (\cite{Ar13}, 2013, communicated by Guicciardini)
and Rabouin--Arthur (\cite{Ra20}, 2020, communicated by Jeremy Gray).
Even earlier, Arthur produced ``Leibniz's Archimedean
infinitesimals''%
\footnote{\label{f5}Among the authors who have pursued such an
Archimedean line is Levey.  Levey is of the opinion that ``Leibniz has
abandoned any ontology of actual infinitesimals and adopted the
\emph{syncategorematic} view of both the infinite and the infinitely
small {\ldots} The interpretation is worth stating in some detail,
both for \emph{propaganda} purposes and for the clarity it lends to
some questions that should be raised concerning Leibniz's
fictionalism'' (\cite[p.\;107]{Le08}; emphasis on ``syncategorematic''
in the original; emphasis on ``propaganda'' added).  In conclusion,
``Leibniz will emerge at key points to be something of an
Archimedean'' (ibid.).  Levey does not elaborate what ``propaganda
purposes'' he may have in mind here, but apparently they don't include
considering Robinson's interpretation as a possibility.  Most
recently, Levey claims that ``Particularly critical for grasping the
conceptual foundation of his method {\ldots} is the fact that for
Leibniz the idea of `infinitely small' is itself understood in terms
of a \emph{variable finite} quantity that can be made
\emph{arbitrarily} small'' (Levey \cite{Le21}, 2021, p.\;142; emphasis
in the original), but provides no evidence.}
(\cite{Ar07}, 2007).

However, Leibniz made it clear that the infinitesimal method is an
\emph{alternative} to, rather than being merely a \emph{fa\c con de
parler} for, Archimedean methods involving exhaustion.  Thus, Leibniz
wrote in his \emph{De Quadratura Arithmetica} (DQA) already in
1675/76:
\begin{quote}
What we have said up to this point about infinite and infinitely small
quantities will appear obscure to certain people, as does everything
new--although we have said nothing that cannot be easily understood by
each of them after a little reflection: indeed, whoever has understood
it will recognize its fecundity.  (Leibniz as translated by
Rabouin--Arthur in \cite[p.\;419]{Ra20})
\end{quote}
How can one overcome such apprehensions concerning the novelty of the
infinite and the infinitely small?  Leibniz offers a way forward by
dispensing with ontological preoccupations as to their existence ``in
nature,'' and treating them as fictional:
\begin{quote}
It does not matter whether there are such quantities in nature, for it
suffices that they be introduced by a fiction, since they provide
abbreviations of speaking and thinking, and thereby of discovery as
well as of demonstration, {\ldots}
(ibid.)
\end{quote}
We deliberately interrupted the passage at a comma for reasons that
will become clear shortly.  Leibniz's rejection of the assumption that
infinitesimals exist in nature occurs in the same sentence as their
description as fictions, indicating that their fictionality is the
contrary of being instantiated in nature.  Leibniz continues:
\begin{quote}
{\ldots} so that it is not always necessary to use inscribed or
circumscribed figures and to infer ad absurdum, and to show that the
error is smaller than any assignable.  (ibid.)
\end{quote}
Leibniz asserts that the infinitesimal method renders it
\emph{unnecessary} to use inscribed or circumscribed figures.  He thus
makes it clear that the infinitesimal method (freed of any ontological
assumptions as to infinitesimals being instantiated in nature)
provides an alternative to, rather than being a \emph{fa\c con de
parler} for -- as per the syncategorematic reading -- Archimedean
methods (involving inscribed and circumscribed figures mentioned in
the passage).  Clearly at odds with the Leibnizian intention in this
passage, the syncategorematic reading specifically makes it
\emph{necessary} to use inscribed and circumscribed figures.%
\footnote{\label{f12}Levey quotes this Leibnizian passage in
\cite[p.\;146]{Le21} and claims that ``Leibniz's view in DQA of the
infinitely small and its correlative idea of the infinitely large as
useful fictions, whose underlying truth is understood in terms of
relations among variable finite quantities that could be spelled out
painstakingly if necessary, will be his longstanding position''
(ibid.), without however providing any evidence.  As we argued, rather
than support the syncategorematic reading, the 1676 passage actually
furnishes evidence against it.}

The fictionalist view as expressed in DQA is consistent with what
Leibniz wrote decades later concerning the interpretation of the
infinitely small.  Thus, in a letter to Varignon, Leibniz commented on
the ideal status of infinitesimals in a way similar to that in DQA:
\begin{quote}
D'o\`u il s'ensuit, que si quelcun%
\footnote{Here and in the following we preserved the spelling as found
in Gerhardt.}
n'admet point des lignes infinies et infiniment petites \`a la rigueur
metaphy\-sique et comme des choses reelles, il peut s'en servir
seurement comme des notions \emph{ideales} qui abregent le
raisonnement, semblables \`a ce qu'on nomme racines imaginaires dans
l'analyse commune{\ldots}%
\footnote{Loemker's English translation: ``It follows from this that
even if someone refuses to admit infinite and infinitesimal lines in a
rigorous metaphysical sense and as real things, he can still use them
with confidence as ideal concepts which shorten his reasoning, similar
to what we call imaginary roots in the ordinary algebra''
\cite[p.\;543]{Le89}.}
(Leibniz to Varignon \cite{Le02}, 2 february 1702, GM \cite{Ge50}, IV,
p.\;92; emphasis added).
\end{quote}
Here Leibniz describes infinitesimals as ideal notions, and expresses
a reluctance to get into metaphysical debates as to their reality.

\subsection
{Was the fictionalist interpretation a late development?}
\label{s14}

Rabouin and Arthur (RA) allege the following:
\begin{enumerate}
\item
Some authors claim that when Leibniz called them `fictions' in
response to the criticisms of the calculus by Rolle and others at the
turn of the century, he had in mind a different meaning of `fiction'
than in his earlier work, involving a commitment to their existence as
non-Archimedean elements of the continuum.  \cite[Abstract]{Ra20}
\item
It has been objected [by such authors] that although Leibniz
characterized infinitesimals as non-existent fictions in DQA and other
writings of the mid-1670s, it cannot be assumed that he continued to
hold the same view of fictions after he had fully developed the
differential calculus.  \cite[pp.\;403--404]{Ra20}
\end{enumerate}
In these passages, RA attribute to their opponents the view that
Leibniz may have had a different meaning of \emph{fiction} when he
fully developed his calculus as compared to his DQA.\, However, RA
provide no evidence to back up the claim that their opponents
``persist in holding that Leibniz's fictionalist interpretation was a
late development, prompted by the criticisms of Rolle and others in
the 1690s.''%
\footnote{In their note~8, RA attribute such a view to Jesseph.
However, Jesseph observes that ``there are certainly traces of
[Leibniz's fictionalist position] as early as the 1670s''
\cite[p.\;195]{Je15}.  RA give an erroneous page for Jesseph's
observation.}

Did Leibniz change his mind about the approach developed in DQA?  It
is true that, already in the late 1670s, Leibniz felt that the new
infinitesimal calculus may have superceded the techniques of DQA.
Thus, in a 1679 letter to Huygens, Leibniz wrote:
\begin{quote}
I have left my manuscript on arithmetical quadratures at Paris so that
it may some day be printed there.%
\footnote{Knobloch quotes only the first sentence ``I have left my
manuscript on arithmetical quadratures at Paris so that it may some
day be printed there'' and infers that ``Leibniz did want to publish
it'' \cite[p.\;282]{Kn17}.  The rest of the quotation makes such an
inference less certain.}
But I have advanced far beyond studies of this kind and believe that
we can get to the bottom of most problems which now seem to lie beyond
our calculation; for example, quadratures, the inverse method of
tangents, the irrational roots of equations, and the arithmetic of
Diophantus.  (Leibniz in 1679 as translated by Loemker
\cite[p.\;248]{Le89}).
\end{quote}
Thus already by 1679, Leibniz ``advanced far beyond'' his DQA.
However, DQA being possibly superceded does not mean that Leibniz
changed his mind about the meaning of infinitesimals as fictions.  And
if Leibniz did change his mind, then RA would be contradicting Leibniz
himself (rather than contradicting their opponents), as they do when
they claim to account for Leibnizian infinitesimals
\begin{quote}
``in keeping with the Archimedean axiom'' \cite[Abstract]{Ra20},
\end{quote}
whereas Leibniz made it clear that his incomparables violate Euclid's
Definition\;V.4 in a pair of 1695 texts \cite{Le95a}, \cite{Le95b}
(see Section~\ref{V4}).

Similarly, Leibnizian bounded infinities are greater than any
assignable quantity (see Section~\ref{s43b} for a modern formalisation).
Leibniz clarified the meaning of the term \emph{incomparable} as
follows:
\begin{quote}
C'est ce qui m'a fait parler autres fois des incomparables, par ce que
ce que j'en dis a lieu \emph{soit} qu'on entende des grandeurs
infiniment petites \emph{ou} qu'on employe des grandeurs d'une
petitesse inconsiderable et suffisante pour faire l'erreur moindre que
celle qui est donn\'ee.  (Leibniz quoted by Pasini \cite{Pa88}, 1988,
p.\;708; emphasis added)
\end{quote}
Leibniz' sentence exploits the structure ``soit {\ldots} ou {\ldots}"
to present a pair of distinct meanings of the term
\emph{incomparable}: (A) as a practical tool enabling one to make an
error smaller than a value given in advance,%
\footnote{\label{f13}The \emph{incomparables} in this sense are
sometimes described by Leibniz as \emph{common}.  Thus, they are
described as \emph{incomparables communs} in the 1702 letter to
Varignon quoted in Section~\ref{s32b}.  The distinction has often been
overlooked by commentators; see note~\ref{f24}.}
and (B) as an infinitesimal properly speaking.%
\footnote{Leibnizian conciliatory methodology (see Mercer \cite{Me06},
  2006) could account for his desire to include a mention of the
  A-method (favored by some contemporaries on traditionalist grounds)
  alongside the B-method.}
The 1695 texts used the term in the sense~(B).\, RA quote this
Leibnizian passage \emph{twice} yet fail to notice that it challenges
their main thesis concerning Leibnizian infinitesimals (see
Section~\ref{s22} for an analysis of the difficulties with RA's
reading).

\subsection{Euclid's Definition\;V\!.4}
\label{V4}

A pair of Leibnizian texts \cite{Le95a}, \cite{Le95b} from 1695
mention Euclid's definition\;{V\!.4} (referred to as V\!.5 by
Leibniz): a letter to l'Hospital and a published response to
Nieuwentijt.  Breger notes the violation of Euclid's definition by
Leibnizian incomparables:
\begin{quote}
In a letter to L'H\^opital of 1695, Leibniz gives an explicit
definition of incomparable magnitudes: two magnitudes are called
incomparable if the one cannot exceed the other by means of
multiplication with an arbitrary (finite) number, and he expressly
points to Definition\;5 of the fifth book of Euclid quoted above.
\cite[p.\;73--74]{Br17}
\end{quote}
See Bair et al.\;(\cite{18a}, 2018, Sections\;3.2--3.4) for a detailed
analysis of these two mentions of V\!.4.

Levey claims that ``in a 1695 response to criticisms of the calculus
by Bernard Nieuwentijt {\ldots} Leibniz explains that his infinitely
small differential (say,~$dx$ in
\mbox{$x+dx$}) is not to be taken to be a fixed small quantity''
\cite[p.\;147]{Le21}.  Levey goes on to translate a Leibnizian passage
as follows:
\begin{quote}
Such an increment cannot be exhibited by construction. Certainly, I
agree with Euclid bk.\;5, defin.\;5, that only those homogeneous
quantities are comparable of which one when multiplied by a number,
that is, a finite number, can exceed the other. And I hold that any
entities whose difference is not such a quantity are equal.  (\ldots)
This is precisely what is meant by saying that the difference is
smaller than any given. (GM 5.322) (Leibniz as translated by Levey in
\cite[p.\;147]{Le21})
\end{quote}
Arthur makes similar claims about this passage in (\cite{Ar13}, 2013,
p.\;562).  However, such claims concerning this Leibnizian passage are
problematic, as signaled in the following seven points.

\medskip
\textbf{1.}  Both Arthur and Levey skip Leibniz's crucial introductory
sentences (quoted in Section~\ref{s28} below) where Leibniz discusses
both his generalized relation of equality up to negligible terms, and
his incomparables.

\medskip
\textbf{2.}  Arthur and Levey overlook the fact that Leibniz is
\emph{not} claiming that homogeneous quantities are always comparable.
Rather, Leibniz merely gives the definition of comparability in terms
of\, V\!.4.  Levey seeks to create a spurious impression that Leibniz
endorses the Archi\-medean property, whereas Leibniz merely gives a
definition of what it means for quantities to be comparable.

\medskip
\textbf{3.}  The promoters of the syncategorematic reading are in
disarray when it comes to interpreting this Leibnizian passage.  RA
\cite[p.\;432]{Ra20} acknowledge the violation of Euclid V\!.4 at
least as a nominal definition of incomparable quantities that may be
used as if they exist under certain specified conditions (see
Section~\ref{s34}).  In a comment appearing in our third epigraph, RA
even summarize the Leibnizian definition of incomparability in modern
notation as follows.  They denote by ``INC'' the relation of being
incomparable, and write:
\begin{quote}
[Leibniz] has defined {\ldots} what it is for one quantity (say,~$dx$)
to be incomparable \emph{in relation to} another (say,~$x$).
Formally, Leibniz's definition gives:\; For a given~$x$
and~$dx$,\,~$dx$ INC~$x$ iff \mbox{$\neg(\exists n) (ndx>x)$}.
(\cite[p.\;433]{Ra20}; emphasis in the original)
\end{quote}
By contrast, Levey \cite[p.\;147]{Le21} and Arthur \cite{Ar13}
interpret the passage directly in an Archime\-dean fashion involving
understanding the clause ``not a comparable difference" as \emph{no
difference} rather than \emph{possibly an incomparable difference}.%
\footnote{As mentioned earlier, this interpretation is untenable since
in the preceding sentences -- not quoted by Levey -- Leibniz defines
his generalized equality where the difference may be incomparably
small rather than necessarily absolutely zero.}

\medskip
\textbf{4.}  The contention that Leibniz had in mind a generalized
relation of equality up to negligible terms is corroborated by the
fact that Jacob Bernoulli's student, Jacob Hermann, gave a similar
definition in his 1700 rebuttal of Nieuwentijt's criticisms:
\begin{quote}
According to Hermann the whole difficulty is based on an ambiguous use
of the two notions `aequalis' and `incomparabilis'.  With regard to
equality, he reiterates Leibniz's definition, formulating it like
this: ``quaecunque data quavis minore differentia differunt, aequalia
esse.'' (``Whatever differs by a difference smaller than any given
quantity is equal'').  (Nagel \cite{Na08}, 2008, p.\;206)
\end{quote}
Hermann also gave a definition of infinitesimals similar to
Leibniz's.%
\footnote{\label{f16}Hermann defines infinitesimals as follows:
``Quantitas vero infinite parva est, quae omni assignabili minor est:
{\&} talis Infinitesima vel Differentiale vocatur''
\cite[p.\;56]{He00}.  Thus an infinitesimal is smaller than every
assignable quantity, just as in Leibniz.}
While the idea of Leibniz being ahead of his time possesses certain a
priori \emph{philosophical} plausibility, \emph{historically} speaking
the idea of Leibniz as a harbinger of the \emph{Epsilontik} must
postulate a problematic discontinuity with both his predecessors
Torricelli and Roberval and his successors l'Hospital, Hermann, the
Bernoullis and others (see also note~\ref{f5b}).
 
\medskip
\textbf{5.}  Levey's opening claim that ``Leibniz explains that his
infinitely small differential {\ldots} is not to be taken to be a
fixed small quantity'' is not supported by the Leibnizian passage at
all, since the passage mentions neither constant nor variable
quantities.

\medskip
\textbf{6.}  In view of Breger's remark quoted at the beginning of the
current Section~\ref{V4}, Levey appears to attribute to Leibniz a
rather paradoxical stance of claiming that his incomparables
violate~{V\!.4} in one 1695 text (letter to l'Hospital), and
satisfy~{V\!.4} in another text (response to Nieuwentijt) from the
same year.  The analysis presented in \cite{18a} indicates that in
both texts, Leibniz made it clear that his incomparables%
\footnote{Not to be confused with \emph{common} incomparables; see
notes~\ref{f13} and \ref{f24}.}
violate~V\!.4.

\medskip
\textbf{7.} Levey follows Ishiguro's narrative in the following sense.
He first presents what he describes as the ``popular history'' of LC,
and then proposes an alternative.  According to the popular history as
described by Levey, Leibniz used unrigorous infinitesimals to advance
geometry without worrying about the foundations of his method.
Levey's alternative account holds that Leibniz was in fact rigorous
(by the standards of his day), but his concept of infinitesimal was
misunderstood: what he really meant was Archimedean variable
quantities.  Thus, Levey writes:
\begin{quote}
Leibniz's calculus was developed with scrupulous care for foundational
matters and was never wedded to the idea of infinitesimals as fixed
quantities greater than zero but less than any finite value.
\cite[p.\;140]{Le21}
\end{quote}
We can agree that LC was developed with scrupulous care for
foundational matters.  However, Levey's alternative account is rooted
in a Weierstrassian concept of rigor,%
\footnote{Additional difficulties with Levey's reading are signaled in
notes \ref{f5}, \ref{f12}, \ref{f26}, \ref{f30} \ref{f44}.}
and overlooks another possibility--namely, that Leibnizian
infinitesimals were rigorous (by the standards of his time), but to
appreciate it one needs to overcome the Weierstrassian limitations and
be open to taking Leibniz at face value when he makes it clear that
his infinitesimals violate Definition V\!.4.  The coherence of LC with
its incomparables and bounded infinities is underscored by the
existence of modern formalisations (see Section~\ref{s24}) faithful to
the Leibnizian procedures.

\subsection
{Questions for dialog}
\label{s15}

While the syncategorematic reading of LC is not the main focus of our
analysis in the present text, we would like to present a few questions
for a possible dialog concerning the Ishiguro--Arthur reading.  In the
recent literature, one detects oblique references of the following
type:
\begin{itemize}
\item
``some recent theories of infinitesimals as non-Archimedean entities''
(Arthur \cite{Ar13}, 2013, p.\;553).
\item
``Advocates of nonstandard analysis routinely refuse to acknowledge
this'' (Spalt \cite{Sp15}, 2015, p.\;121).
\item
``Recently there have been attempts to argue that Leibniz, Euler, and
even Cauchy could have been thinking in some informal version of
rigorous modern non-standard analysis, in which infinite and
infinitesimal quantities do exist.  However, a historical
interpretation such as the one sketched above that aims to understand
Leibniz on his own terms, and that confers upon him both insight and
consistency, has a lot to recommend it over an interpretation that has
only been possible to defend in the last few decades''%
\footnote{Apart from its failure to either cite or name scholars being
criticized, Gray's position endorsing historical authenticity seems
reasonable--until one examines Gray's own historical work, where
Euler's foundations are claimed to be ``dreadfully weak''
\cite[p.\;6]{Gr08} and ``cannot be said to be more than a gesture''
\cite[p.\;3]{Gr15}, whereas Cauchy's continuity is claimed to be among
concepts Cauchy defined using ``limiting arguments''
\cite[p.\;62]{Gr08b}.  Such claims are arguably symptomatic of
butterfly-model thinking (see Section~\ref{s17b}).  For a more
balanced approach to Euler and Cauchy see respectively \cite{17b} and
\cite{20a}.}
(Gray \cite{Gr15}, 2015, p.\;11).
\item
``Certain scholars of the calculus have denied that the interpretation
of infinitesimals as syncategorematic was Leibniz's mature view, and
have seen them as fictions in a different sense''
(Arthur\;\cite{Ar18}, 2018, p.\;156).
\item
``on se gardera n\'eanmoins de forcer un parall\`ele que Leibniz ne
place pas l\`a o\`u l'interpr\'etation formaliste le place'' (Rabouin
in \cite{Le18}, 2018, p.\;96).
\end{itemize}
Alas, such references come unequipped with either quotation or
citation.  These authors apparently feel that the answer to Probst's
question (see Section~\ref{one}) is evident.  It is to be hoped that a
more meaningful dialog with proponents of the Ishiguro--Arthur reading
and its variants can be engaged in the future.%
\footnote{It is commendable that in a recent text, Rabouin and Arthur
(\cite{Ra20},\;2020) engage their opponents in a more open fashion
than has been the case until now, even though we raise many questions
concerning the quality of their engagement in Sections~\ref{s14}
and~\ref{s22}.}
Such a dialog could use the following questions as a starting point:
\begin{enumerate}
\item
\emph{Were there two separate methods in the Leibnizian calculus or
only one?}  Bos claims that there were two (one involving exhaustion
in an Archi\-me\-dean context, and a separate one relying on the law
of continuity and exploiting inassignable infinitesimals), and
provides textual evidence.  A number of authors disagree with Bos,
including Arthur \cite[p.\;561]{Ar13}, Breger
\cite[pp.\;196--197]{Br08}, and Spalt \cite[p.\;118]{Sp15}.%
\footnote{Rabouin and Arthur (\cite{Ra20}, 2020) seem to acknowledge
the presence of two methods, but then go on to claim that method B
exploiting infinitesimals~$dx$ and~$dy$ is easily transcribed in terms
of assignable quantities~$(d)x$ and~$(d)y$.  The mathematical
coherence (or otherwise) of such a claim is analyzed in
Section~\ref{s42}.}
\item
\emph{Are there hidden quantifiers lurking behind Leibnizian
calculus?}  Ishiguro said there were (see Section~\ref{s13}).  Rabouin
seeks to distance himself from hidden quantifiers in
\cite[p.\;362]{Ra15} but goes on to endorse Ishiguro in note 25 on the
same page.  Breger endorses ``absorbed'' quantifiers in
\cite[p.\;194]{Br08}.
\item
\emph{Is it or is it not appropriate to use Robinson's framework in
order to formalize Leibnizian analysis?}  Here there seems to be a
uniform agreement among scholars including 
Arthur, Bassler, Breger, Gray, Ishiguro, Levey, Rabouin \cite[note\;2,
pp.\;95--96]{Le18}, and Spalt
that a Robinsonian interpretation is not possible.  Yet Arthur
(\cite{Ar13}, 2013) compares LC to Smooth Infinitesimal Analysis where
infinitesimals similarly exhibit non-Archimedean behavior.  Granted he
finds some differences in addition to similarities, but he does not
seem to accord a status of ``comparison is possible'' to Robinson's
framework.  The reasons for this need to be understood.
\item
\emph{Do scholars following Ishiguro tend to apply a modern framework
to interpreting Leibniz, namely the Weierstrassian one?}  In some of
his later writings, Knobloch explicitly describes Leibnizian
infinitesimals in terms of Weierstrass and the \emph{Epsilontik}; see
e.g., \cite[pp.\;2,\;7]{Kn21}.
\item
\emph{Must Leibnizian infinitesimals be interpreted as limits?}
Guicciardini (\cite{Gu18}, 2018, p.\;90) replies in the affirmative.
\end{enumerate}
Philosophical \emph{partis pris} among historians of mathematics have
been analyzed by Hacking; see Section~\ref{s17b}.

\subsection
{Hacking's butterfly/Latin models; Carnap's `empty words'}
\label{s17b}

Some historians may be conditioned by their undergraduate mathematical
training to interpret Leibniz in anachronistic ways.  Thus, in one of
his later articles, Knobloch writes:
\begin{quote}
The set of all finite cardinal numbers~$1, 2, 3, \ldots$ is a
transfinite set.  Its cardinal number is Alef$_0$.  This is the least
cardinal number being larger than any finite cardinal number.
Leibniz's terminology implies actual infinity though he rejects the
existence of an infinite number, etc.  (Knobloch \cite{Kn18}, 2018,
pp.\;13--14)
\end{quote}
The assumption that the only way to formalize infinity is by means of
Cantorian infinite cardinalities creates tensions with what Leibniz
wrote.  Specifically, \emph{infinita terminata} can be understood
without Cantorian infinite cardinalities; see Section~\ref{s43}.

Some Leibniz scholars work under the related assumption that
`Archi\-medean' is the only natural way of thinking of the continuum.
Thus Arthur and Rabouin endeavor to clear Leibniz of any suspicion of
entertaining unnatural thoughts such as infinitesimals.  They would
rather declare Leibniz to be mired in contradictions (see Section
\ref{s31b}) than admit that he dealt with non-Archimedean phenomena.

Such ``conceptually repugnant'' attitudes are illustrated in our
second epigraph taken from Levey \cite[p.\;55, note\;9]{Le98}.  Such
scholars fail to take into account the contingency of the historical
development of mathematics emphasized by Hacking (\cite{Ha14}, 2014,
pp.\;72--75, 119).

Hacking contrasts a model of a deterministic (genetically determined)
biological development of animals like butterflies (the
egg--larva--cocoon--butterfly sequence), with a model of a contingent
historical evolution of languages like Latin.  Emphasizing determinism
over contingency in the historical evolution of mathematics can easily
lead to anachronism.%
\footnote{For a case study in Cauchy scholarship see Bair et
al.\;(\cite{19a}, 2019).  Opposition to Marburg neo-Kantians' interest
in infinitesimals is documented in (\cite{13h}, 2013).}
With reference to infinitesimal calculus, the Latin model challenges
the assumption that the so-called search for rigor inevitably led to
the arithmetisation of analysis, accompanied by the banning of
infinitesimals and the establishment of Weierstrassian foundations
ultimately formalized in Zermelo--Fraenkel set theory (ZF) in
the~$\in$-language (formalizing the membership relation).  The
development towards the Weierstrassian limit concept seemed inevitable
to Carnap who wrote:
\begin{quote}
Leibniz and Newton {\ldots} thought that they had a definition which
allowed them to have a conceptual understanding of ``derivative".
However, their formulations for this definition used such expressions
as ``infinitesimally small magnitudes'' and quotions [sic] of such,
which, upon more precise analysis, turn out to be pseudo concepts
(empty words).  It took more than a century before an unobjectionable
definition of the general concept of a limit and thus of a derivative
was given.  Only then all those mathematical results which had long
since been used in mathematics were given their actual meaning.
\cite[p.\;306--307]{Ca03}
\end{quote}
Carnap's assumption is that the ``actual meaning'' of a mathematical
concept like the derivative could not possibly involve ``empty words''
like infinitesimals.  In 2020 one finds a historian who claims that
``Euler was not entirely successful in achieving his aim since he
introduced infinitesimal considerations in various proofs''
\cite[p.\;11]{Fe20}.  Such assumptions obscure the possibility that
the so-called arithmetisation of analysis could have occurred while
retaining infinitesimals.  Hacking's Latin model enables scholars to
contemplate alternative possibilities of historical development of
mathematics.  \emph{Pace} Carnap, set theories in the more versatile
$\st$-$\in$-language (formalizing both the Leibnizian
assignable/inassignable distinction and the membership relation)
provide more faithful conceptualisations of the procedures of LC than
the theory ZF; see Section~\ref{s24} for further details.

\subsection{Grafting of {\law} on the calculus of Leibniz}
\label{s17}

Richard Arthur compared LC to a modern theory of infinitesimals,
namely Smooth Infinitesimal Analysis ({\law}).  Arthur claims to find
``many points in common'' (in addition to differences) between LC and
{\law}:
\begin{quote}
I then turn to a comparison of Leibniz's approach with the recent
theory of infinitesimals championed by John Bell, Smooth Infinitesimal
Analysis ({\law}), of which I give a brief synopsis in Sect.\;3.  As we
shall see, this has many points in common with Leibniz's approach: the
non-punctiform nature of infinitesimals, their acting as parts of the
continuum, the dependence on variables (as opposed to the static
quantities of both Standard and Non-standard Analysis), the resolution
of curves into \emph{infinite-sided polygons}, and the finessing of a
commitment to the existence of infinitesimals.%
\footnote{Arthur's finessing comment is analyzed in
Section~\ref{s31}.}
(Arthur \cite{Ar13}, 2013, p.\;554; emphasis added)
\end{quote}
The claimed similarities including the following issues:
\begin{enumerate}
[label={~(AR\theenumi)}]
\item
\label{ra1}
\hskip36pt non-punctiform nature of the continuum;
\item
\label{ra2}
\hskip36pt
variables;
\item
\label{ra3}
\hskip36pt
infi\-nite-sided polygons;
\item
\label{ra4}
\hskip36pt
fictionality of infinitesimals.
\end{enumerate}
We will argue that issue \ref{ra1} is largely irrelevant to
interpreting the procedures of LC, while Arthur's case for \ref{ra2}
is based on equivocation.  Furthermore, Robinson's framework is more
successful than {\law} on issues \ref{ra3} and \ref{ra4}.

We present a detailed analysis of Arthur's comparison of LC to {\law}
in Section~\ref{s2}.  The position of Rabouin and Arthur is analyzed
in Section~\ref{s22}.  We propose a formalisation of the procedures of
the Leibnizian calculus in modern mathematics in Section~\ref{s24}.
Some technical details on {\law} are reviewed in Section~\ref{s3}.  In
Section~\ref{s6} we focus on interpretations of the seminal text
\emph{Cum Prodiisset}.

\section{Leibniz and Smooth Infinitesimal Analysis}
\label{s2}

Arthur developed a comparison of LC with {\law} in (\cite{Ar13},
2013).  In this section we will analyze Arthur's comparison.  Any
discussion of {\law} should mention the foundational work of Lawvere
starting in 1967 (see e.g., \cite{La80}, 1980) and the books by
Moerdijk and Reyes (\cite{Mo91}, 1991) and by Kock (\cite{Ko06},
2006), among other authors.

\subsection{Knobloch on Arthur and {\law}}
\label{s21}

Knobloch analyzed Arthur's comparison of LC and {\law} in a
\emph{Zentralblatt} review:
\begin{quote}
[Richard Arthur] compares Leibniz's infinitesimals with those of
Bell's smooth infinitesimal analysis ({\law}) published in
1998. Hence, in Section\;3, he gives a brief synopsis of {\law}.\,
Both notions of infinitesimals use non-punctiform infinitesimals and
the \emph{resolution of curves into infinite-sided polygons}.  The
Leibnizian polygonal representation of curves is closely related to
Bell's principle of microstraightness.  Yet, there are crucial
differences.  (Knobloch \cite{Kn13}, 2013; emphasis added)
\end{quote}
Note that {\law} depends on a category-theoretic foundational
framework and on intuitionistic logic to enable nilpotency of
infinitesimals (see Section~\ref{s31}).  Thus Arthur seeks to rely on
the resources of modern mathematics, including non-classical logics,
in a comparison with LC.\, While we welcome such a display of
pluralism on Arthur's part, we also agree with both John L. Bell and
Arthur that {\law} is closer to Nieuwentijt's approach to calculus
than to Leibniz's (see further in Section~\ref{s25}).  All {\law}
nilpotent infinitesimals~$\epsilon$ have the property that
both~$\epsilon$ and~$-\epsilon$ are smaller than~$\frac{1}{n}$:
\[
\text{for all\;} \epsilon \text{\;in\;} \Delta,%
\footnote{The definition of the part~$\Delta$ appears in
Section~\ref{s33}.}
\text{\ one has\ } \epsilon <\frac1n,
\]
signifying non-Archimedean behavior; see (Bell \cite{Be08}, 2008,
Exercise\;8.2, p.\;110).

Leibniz endorsed what is known today as the law of excluded middle in
a passage written around 1680 in the following terms:

\begin{quote}
Every judgment is either true or false. No judgment is simultaneously
true and false.  Either the affirmation or the negation is true.
Either the affirmation or the negation is false.  For every truth a
reason can be provided, excepting those first truths in which the same
thing is affirmed of the thing itself or is denied of its
opposite. ~$A$ is~$A$.  \mbox{$A$ is not not-$A$}.  (Leibniz as
translated by Arthur in \cite[pp.\;237--239]{Ar01})
\end{quote}
The assumption that proposition~$A$ is equivalent to~$\neg\neg A$ is
the law of excluded middle, rejected by intuitionists.  Arthur
translated and edited this Leibnizian passage in 2001 but did not
report it in his text (\cite{Ar13}, 2013) comparing LC to the theory
SIA based on intuitionistic logic.  The following question therefore
arises.

\begin{quote}
\emph{Question.}  Why does Arthur's voluminous output on the
Leibnizian calculus systematically eschew readings based on Robinson's
framework \cite{Ro66} or the foundational approaches of Hrbacek
\cite{Hr78} and Nelson \cite{Ne77},
%
%
%
based as they are on classical logic, and enabling straightforward
transcription (see Section~\ref{s24}) of both Leibniz's
assignable/inassignable dichotomy and his infinitely many orders of
infinitesimal and infinite numbers?
\end{quote}
An attempt by Rabouin and Arthur to address this question is analyzed
in Section~\ref{s22}.

Meanwhile, Arthur appears to avoid Robinsonian infinitesimals as
zealously as Leibniz avoided atoms and material indivisibles (for an
analysis of Leibniz on atoms and indivisibles see Bair et al.\
\cite{18a}, 2018, Section\;2.6).

\subsection{\emph{Infinita terminata}}
\label{s23}

A key Leibnizian distinction between bounded and unbounded infinity is
sometimes insufficiently appreciated by commentators (including
Rabouin and Arthur).  In DQA, Leibniz contrasts bounded infinity and
unbounded infinity in the following passages:
\begin{quote}
But as far as the activity of the mind with which we measure infinite
areas is concerned, it contains nothing unusual because it is based on
a certain fiction and proceeds effortlessly on the assumption of a
certain, though \emph{bounded, but infinite line}; therefore it has no
greater difficulty than if we were to measure an area that is finite
in length.  (Leibniz, DQA \cite{Le04b}, Scholium following
Propositio\;XI; translation ours; emphasis added)
\end{quote}
Leibniz elaborates as follows:
\begin{quote}
Just as points, even of infinite numbers, are unsuccessfully added to
and subtracted from a bounded line, so a bounded line can neither form
nor exhaust an unbounded one, however many times it has been repeated.
This is different with a \emph{bounded but infinite line} thought to
be created by some multitude of finite lines, although this multitude
exceeds any number.  And just as a \emph{bounded infinite line} is
made up of finite ones, so a finite line is made up of infinitely
small ones, yet divisible.%
\footnote{It is worth noting that bounded/unbounded is not the same
distinction as potential/actual infinity.  Bounded infinity is a term
Leibniz reserves mainly to characterize the quantities used in his
infinitesimal calculus, namely the infinitely small and infinitely
large.}
(ibid.)
\end{quote}
The importance of the distinction is stressed by Mancosu (\cite{Ma96},
1996, p.\;144).  RA appear to be aware of the distinction and even
mention the term \emph{terminata} three times in their article, but
don't fully appreciate its significance.%
\footnote{Furthermore, RA misinterpret Leibnizian bounded infinities
when they compare them to compactifications in modern mathematics; see
Section~\ref{s38}.}
We propose a possible formalisation of the bounded/unbounded
distiction in modern mathematics in Section~\ref{s43}.

\subsection{Infinite-sided polygons}
\label{s25}

Arthur claims that a common point between LC and {\law} consists in
viewing a curve as a polygon with infinitely many sides (see
Section~\ref{s17}).  Knobloch takes note of such alleged similarity in
his review of Arthur quoted in Section~\ref{s21}.  Are infinite-sided
polygons a point in common between LC and {\law}, as Arthur claims?  A
comparison of LC and {\law} based on infinite-sided polygons would
apparently require both sides of the comparison to envision such
polygons.  Arthur goes on to mention
\begin{quote}
the quotation Bertoloni Meli gives from Leibniz's letter to Claude
Perrault {\ldots} written in 1676: `I take it as certain that
everything moving along a curved line endeavours to escape along the
tangent of this curve; the true cause of this is that \emph{curves are
polygons with an infinite number of sides}, and these sides are
portions of the tangents\ldots' (quoted from Bertoloni Meli
1993,\;75).  (Arthur \cite{Ar13}, 2013, note\;17, pp.\;568--569;
emphasis added).
\end{quote}
Thus far, we have sourced Leibnizian infinite-sided polygons.  What
about the {\law} side of Arthur's comparison?  Arthur's source Bell
does mention \emph{infinilateral polygons} in his article ``Continuity
and infinitesimals" (Bell \cite{Be05}, 2005--2013).  However, Bell
attributes such polygons to Leibniz and Nieuwentijt rather than to
{\law}.  We give three examples:
\begin{enumerate}
\item
``The idea of considering a curve as an infinilateral polygon was
employed by a number of \ldots\;thinkers, for instance, Kepler,
Galileo and Leibniz.''  (Bell \cite{Be05})
\item
``Leibniz's definition of tangent employs both infinitely small
distances and the conception of a curve as an infinilateral polygon.''
(ibid.)
\item
``Nieuwentijdt's infinitesimals have the property that the product of
any pair of them vanishes; in particular each infinitesimal is
`nilsquare' in that its square and all higher powers are zero. This
fact enables Nieuwentijdt to show that, for any curve given by an
algebraic equation, the hypotenuse of the differential triangle
generated by an infinitesimal abscissal increment~$e$ coincides with
the segment of the curve between~$x$ and~$x + e$.  That is, a curve
truly \emph{is} an infinilateral polygon.''  (ibid.; emphasis in the
original)
\end{enumerate}

While Bell speaks of Leibniz and Nieuwentijt as using infinilateral
polygons, Bell's article is silent on infinilateral polygons in
connection with {\law}.  

Similarly, Bell's book on {\law} (\cite{Be08}, 2008) makes no mention
of infinilateral polygons.  There is a good reason for Bell's silence.
Nilsquare infinitesimals, say~$\epsilon$, are by definition not
invertible.  Therefore in {\law} one cannot view a curve as made up
of, say,~$\frac{1}{\epsilon}$ sides or something of that order.
Rather, in {\law}, a curve microlocally ``looks like'' a straight
infinitesimal segment in the sense that, microlocally, a function
coincides with its first Taylor polynomial; see formula~\eqref{e32} in
Section~\ref{s33}.  Therefore infinite-sided polygons can hardly be
said to be a point in common between LC and {\law}.  Arthur's claim to
the contrary is dubious.

\subsection{Bell's microquantities}
\label{s26}

In his more recent book (\cite{Be19}, 2019) published six years after
Arthur's article \cite{Ar13}, Bell mentions infinilateral polygons in
the context of his discussion of Stevin, Kepler, Galileo, Barrow,
Leibniz, l'H\^opital, and Nieuwentijt in his historical Chapter\;2.

In Chapter\;10, Bell gives his definition of \emph{microstraightness}
in terms of the set (part) of microquantities~$\Delta$%
\footnote{The definition of the part~$\Delta$ appears in
Section~\ref{s33} below.}
as follows:
\begin{quote}
If we think of a function~$y = f(x)$ as defining a curve, then, for
any~$a$, the image under~$f$ of the ``microinterval''~$\Delta+a$
obtained by translating~$\Delta$ to~$a$ is straight and coincides with
the tangent to the curve at~$x = a$ {\ldots} In this sense {\ldots}
each curve is ``microstraight''.  \cite[p.\;241]{Be19} (emphasis in
the original)
\end{quote}
Bell goes on to claim that
\begin{quote}
closed curves can be treated as infinilateral polygons, as they were
by Galileo and Leibniz (op. cit., note\;18, p.\;241)
\end{quote}
However, he provides no further elaboration.  The implied assumption
that Leibnizian infinilateral polygons were akin to {\law}'s seems
dubious; to justify such an assumption one would have to stretch the
meaning of the term ``infinilateral'' beyond its apparent meaning of
``possessing infinitely many sides.''  This is because the
\emph{microquantities} in question are not invertible, as already
mentioned in Section~\ref{s25}.  The matter is analyzed further in
Section~\ref{s27}.

\subsection{Infinilateral or infinitangular?}
\label{s27}

In his seminal 1684 text \emph{Nova Methodus}, Leibniz refers to a
curve with infinitely many \emph{angles}, rather than infinitely many
sides.  

While curves in {\law} are microstraight, enabling many elegant
calculations, they don't fit the Leibnizian description.  Having an
angle requires having a vertex, but from the {\law} viewpoint, curves
not only do not have infinitely many vertices but they have no
vertices (where one could speak of an \emph{angle} not\;not different
from~$\pi$) at all.  Perhaps it is for this reason that Bell can only
speak of infini\emph{lateral} polygons rather than
infinit\emph{angular} polygons.  How does Bell handle the issue?  Bell
writes:

\begin{quote}
Leibniz's definition of tangent employs both infinitely small
distances and the conception of a curve as an \emph{infinilateral}
polygon:
\begin{quote}
We have to keep in mind that to find a tangent means to draw a line
that connects two points of a curve at an infinitely small distance,
or the continued side of a polygon with an \emph{infinite number of
angles}, which for us takes the place of a curve.
\end{quote}
In thinking of a curve as an \emph{infinilateral} polygon {\ldots} the
abscissae~$x_0, x_1, \ldots$ and the ordinates~$y_0, y_1$, are to be
regarded as lying infinitesimally close to one another; {\ldots} (Bell
\cite{Be19}, 2019, p.\;70; emphasis added)
\end{quote}
The indented Leibnizian passage is from (Leibniz \cite{Le84}, 1684,
p.\;223).
%
%
%
%
%
%
%
%
Leibniz speaks of infinit\emph{angular} polygons, whereas Bell speaks
of infini\emph{lateral} polygons both before and following the
Leibnizian passage, without commenting on the discrepancy.  Leibnizian
infinitangular polygons are arguably modeled less faithfully in {\law}
than in Robinson's framework, as we argue in Section~\ref{s28}.

\subsection
{Infinitangular and infinilateral polygons in Robinson}
\label{s28}

We saw in Section~\ref{s26} that infinite-sided polygons are not a
feature of the {\law} framework in the context of nilsquare
microquantities.  Meanwhile, in Robinson's framework, an
infinite-sided polygon%
\footnote{More precisely, an internal polygon with an infinite
hyperinteger number of vertices.  In Nelson's framework or in the
theory SPOT one considers a polygon with~$\mu$ sides for a
nonstandard~$\mu\in\N$; see Section~\ref{s46}.}
can easily be chosen to approximate a given curve.  This enables one
to determine the usual geometric entities, such as the tangent line
from a pair of adjacent vertices or the curvature from a triple of
adjacent vertices of the polygon.

To provide a modern interpretation, consider the problem of
determining the tangent line to a curve.  Take a pair of infinitely
close distinct points~$D,E$ on the curve, and consider the line~$DE$.
One can normalize the equation of the line~$DE$ as~$ax+by=c$
where~$a^2+b^2=1$.  Then one can take the standard part (see
Section~\ref{s35}) of the coefficients of the equation of the
line~$DE$ to obtain the equation~$a_o x+b_o y=c_o$ of a true tangent
line to the curve at the standard point~$D_o=E_o$, where~$r_o$ is the
standard part of~$r$ for each of~$r=a,b,c,D,E$.  In this sense, the
line~$DE$ is an approximation to the true tangent line, meaning that
they coincide \emph{up to negligible terms}.  For more advanced
applications see e.g., Albeverio et al.\;(\cite{Al86}, 1986).

Approximation procedures are not unusual for Leibniz.  Leibniz often
mentioned his use of a generalized notion of equality.  Thus, in his
response to Nieuwentijt Leibniz wrote:
\begin{quote}
Furthermore I think that not only those things are equal whose
difference is absolutely zero, but also those whose difference is
incomparably small.  And although this [difference] need not
absolutely be called Nothing, neither is it a quantity comparable to
those whose difference it is.  (Leibniz \cite{Le95b}, 1695, p.\;322)
\end{quote}
Bell makes the following claims: 
\begin{quote}
In {\law}, [1] all curves are microstraight, and [2] closed curves
[are] infinilateral polygons.  [3] Nothing resembling this is present
in NSA [Nonstandard Analysis].  (Bell \cite{Be19}, 2019, p.\;259;
numerals [1], [2], [3] added)
\end{quote}
Bell's claim [1] is beyond dispute.  However, both of Bell's
claims\;[2] and~[3] are dubious.  Claim\;[2] was already analyzed in
Section~\ref{s26}.  Claim\;[3] may be literally true in the sense that
in Robinson's framework, the infinite-sided polygon does not literally
coincide with the closed curve but is only an approximation to the
curve.  However, the approximation is good enough to compute all the
usual entities such as tangent line and curvature, as noted above.
%
%
Leibniz defined the curvature via (the reciprocal of) the radius of
the osculating circle (see Bos \cite{Bo74}, 1974, p.\;36).  The
osculating circle can be defined via the circle passing through a
triple of infinitely close points on the curve (alternatively, the
center of the osculating circle to the curve at a point~$p$ can be
defined via the intersection of two ``consecutive'' normals to the
curve near~$p$).  Johann Bernoulli's approach to osculating circles is
analyzed by Bl{\aa}sj\"o (\cite{Bl17}, 2017, pp.\;89--93).  The
calculation of curvature from a triple of infinitely close points on
the curve is straightforward in Robinson's framework (or its axiomatic
formulations; see Section~\ref{s24}).
%
%


\subsection
{Leibniz to Wallis on Archimedes}

In connection with Leibniz's reference to Archimedes in a letter to
Wallis, Arthur comments as follows:
\begin{quote}
The strict proof operating only with assignable quantities justifies
proceeding by simply appealing to the fact that~$dv$ is incomparable
with respect to~$v$: in keeping with the \emph{Archimedean axiom}, it
can be made so small as to render any error in neglecting it smaller
than any given.
(Arthur \cite{Ar13}, 2013, pp.\;567--568; emphasis added)
\end{quote}
Here Arthur is assuming the~$dv$ is a \emph{common} incomparable and
therefore assignable.%
\footnote{See Section~\ref{s14} and note~\ref{f24} on the distinction
between incomparables and common incomparables.}
Is the Leibnizian~$dv$ an assignable quantity, as Arthur claims?
Arthur goes on to quote Leibniz as follows:
\begin{quote}
Thus, in a letter to Wallis in 1699, Leibniz justifies the rule
for~$d(xy) = xdy + ydx$ as follows[:] \vskip6pt
\begin{quote}
\ldots there remains~$xdy + ydx + dxdy$.  But this [term]~$dxdy$
should be rejected, as it is incomparably smaller than
\mbox{$xdy+ydx$}, and this becomes~$d(xy) = xdy + ydx$, inasmuch as,
if someone wished to \emph{translate} the calculation into the style
of \emph{Archimedes}, it is evident that, when the thing is done using
assignable qualities, the error that could accrue from this would
always be smaller that any given.
\end{quote}
(ibid.; emphasis on ``translate'' added)
\end{quote}
%
%
%
When Leibniz speaks of the possibility of using assignable quantities
only, he explicitly refers to the possibility of a \emph{translation}
of the calculation into the style of Archimedes, rather than to the
original infinitesimal calculation itself.%
\footnote{\label{f35}The Leibnizian passage is also quoted by RA
\cite[p.\;437]{Ra20} who similarly fail to account for the fact that
talking about a \emph{translation} into the style of Archimedes
entails the existence of a separate method exploiting infinitesimals
\emph{\`a la rigueur}.}
Arthur's inference of alleged Archimedean nature of Leibnizian~$dv$ is
dubious.%
\footnote{Furthermore, Arthur's attempt to account for the Leibnizian
derivation of the law~$d(xy)=xdy+y dx$ in Archimedean terms, by means
of quantified variables representing~$dx$,~$dy$,~$d(xy)$ taking
assignable values and tending to zero, would run into technical
difficulties.  Since all three tend to zero, the statement to the
effect that ``the error that could accrue from this would always be
smaller than any given'' understood literally is true but vacuous:
$0=0+0$.  To assign non-vacuous meaning to such an Archimedean
translation of Leibniz, one would have to rewrite the formula, for
instance, in the form~$\frac{d(xy)}{dx}=x\frac{dy}{dx}+ y$.  While
such a paraphrase works in the relatively simple case of the product
rule, it becomes more problematic for calculations involving, e.g.,
transcendental functions.}
%
%

In a letter to Wallis two years earlier, Leibniz already emphasized
the distinction between Archimedean and non-Archimedean techniques:
\begin{quote}
Leibniz's response [to Wallis] was first to distinguish between two
kinds of tetragonistic methods, both of which could be traced back to
Archi\-medes, and then to distinguish these from recent infinitesimal
techniques. Of the older established methods one, he writes, considers
geometrical figures and bodies to be collections of an infinite number
of quantities each of which is incomparably smaller than the whole,
while in the other,%
\footnote{We added a comma for clarity.}
quantities remain comparable to the whole and are taken successively
in infinite number so as eventually to exhaust that whole.  {\dots}
Effectively, the concept of analysis with \emph{the infinite as its
object} became the distinguishing factor between the old and the new.
(Beeley \cite{Be13}, 2013, p.\;57; emphasis added)
\end{quote}

\section{Rabouin--Arthur indictment of infinitesimals}
\label{s22}

Rabouin and Arthur (RA) address the \emph{Question} posed at the end
of Section~\ref{s21} by purporting to identify a factor that
``distinguishes Leibniz's position from Robinson's'' (\cite{Ra20},
2020, p.\;407):
\begin{quote}
For what Robinson demonstrated with his non-standard analysis was
precisely that one can introduce infinite numbers and infinitesimals
\emph{without contradiction}.  The whole point of any `non-standard'
approach is the building of a non-standard model in which what is said
about entities such as `infinite numbers' or `infinitely small
quantities' has to be literally and rigorously true.  By contrast,
Leibniz repeatedly claims that what is said about infinitesimals is
\emph{not true} `\`a la rigueur' and entails a paradoxical way of
speaking.$^{18}$\, (ibid.; emphasis added)
\end{quote}
(More on RA's footnote 18 in Section~\ref{s33c} below.)  Here RA seek
to identify a discrepancy between the approaches of Leibniz and
Robinson.  They claim to find one in terms of the presence of a
contradiction (in the notion of infinitesimal) in Leibniz but not in
Robinson.  Such an indictment of Leibnizian infinitesimals entails the
claim that allegedly only Robinson ``can introduce infinite numbers
and infinitesimals without contradiction'' according to RA.\, Note,
however, the following eight items.

\subsection{Do infinitesimals contradict the part-whole axiom?}
\label{s31b}

RA's underlying assumption is that Leibniz views infinitesimals as
being contradictory, just as unbounded infinite collections (taken as
a whole) are.  In fact, RA's purported discrepancy between Leibnizian
and Robinsonian infinitesimals depends crucially on such an assumption
(Leibnizian infinitesimals are contradictory; Robinson's
aren't;\;QED).\, However, what Leibniz repeatedly described as
contradictory is only \emph{unbounded infinity} seen as a whole (and
contradicting the part-whole axiom), not \emph{bounded infinity}
(\emph{infinita terminata}, such as the reciprocal of an
infinitesimal), which he treated as a fiction because it does not
exist in nature (see Section~\ref{s23}).  Since the issue is not the
existence or otherwise of Leibnizian infinitesimals but rather their
being contradictory or not, RA have failed to show that Leibniz viewed
infinitesimals as contradicting the part-whole axiom.

\subsection{Did Leibniz envision infinity \emph{\`a la rigueur}?}
\label{s32b}

RA allege that Leibniz claimed that ``what is said about
infinitesimals is not true `\`a la rigueur.'\,'' However, Leibniz did
allow for the possibility of infinities \emph{\`a la rigueur} as for
instance in the following passage: ``Et c'est pour cet effect que
j'ay%
\footnote{Here and in the following we preserved the original
spelling.}
donn\'e un jour des lemmes des incomparables dans les Actes de
Leipzic, qu'on peut entendre comme on vent [i.e., `veut'], soit des
infinis \emph{\`a\;la\;rigueur}, soit des grandeurs seulement, qui
n'entrent point en ligne de compte les unes au prix des autres.  Mais
il faut considerer en m\^eme temps, que ces incomparables communs%
\footnote{\label{f24}Such \emph{common} incomparables are not to be
confused with the inassignable ones; see Section~\ref{s14}.  Loemker
gives the following translation of this sentence: ``[T]hese
incomparable magnitudes themselves, as commonly understood, are not at
all fixed or determined but can be taken to be as small as we wish in
our geometrical reasoning and so have the effect of the infinitely
small in the rigorous sense'' \cite[p.\;543]{Le89}.  The translation
is not successful as it obscures the fact that \emph{communs} modifies
\emph{incomparables}.  Leibniz is not speaking of common understanding
of incomparable quantities, but rather of common incomparables as
opposed to bona fide ones.  Horv\'ath in \cite[p.\;66]{Ho86} provides
both the original French and an English translation, but omits the
adjective \emph{communs} in his translation.}
m\^emes n'estant nullement fixes ou determin\'es, et pouvant estre
pris aussi petits qu'on veut dans nos raisonnemens Geometriques, font
l'effect des infiniment petits rigoureux'' (Leibniz \cite{Le02}, 1702,
p.\;92; emphasis added).  Leibniz' ``soit {\ldots} soit {\ldots}''
construction clearly indicates that infinity \emph{\`a la rigueur} is
one of the possibilities.  The precise meaning Leibniz attached to
such \emph{infinis \`a la rigueur} is fictional quantities (lacking
instantiation in nature) greater than any assignables (see
Section~\ref{s14}).

\subsection{Leibniz to Wolff on fictions}
\label{s33c}

Concerning RA's footnote 18 (attached to the passage reproduced at the
beginning of the current Section~\ref{s22}), we note the following.
Here RA seek to source their claim that when Leibniz referred to
infinitesimals as fictional, he meant that they involve a
contradiction.  To this end, they quote a lengthy passage from a 1713
letter from Leibniz to Wolff.  Here Leibniz writes: ``[E]ven though in
my opinion [the infinitely small] encompass something of the
\emph{fictive and imaginary}, this can nevertheless be rectified by a
reduction to ordinary expressions so readily that no error can
intervene'' (GM\;V\;385) (Leibniz as translated by RA in \cite{Ra20},
note\;18).  In the passage quoted, Leibniz makes no mention of
contradictions, and merely describes infinitesimals as \emph{fictive
and imaginary.}  If RA wish to demonstrate by this that Leibniz viewed
\emph{fictions} (such as infinitesimals) as contradictory then their
demonstration is circular.  Thus, the passage does not substantiate
RA's claim.

\subsection{What kind of definition?}
\label{s34}

RA claim that the Leibnizian definition of incomparables is
\emph{nominal} rather than \emph{semantic}.  While they acknowledge
that it is ``clear and distinct'' \cite[p.\;406]{Ra20}, they claim
that it does not assert existence or even possible existence.
According to RA, Leibniz' theory of knowledge enables one to derive
truths using contradictory concepts,%
\footnote{To buttress such a claim, RA offer what they describe as ``a
short digression on Leibniz's theory of knowledge''
\cite[p.\;406]{Ra20}.  The evidence they provide is the following
Leibnizian passage: ``For we often understand the individual words in
one way or another, or remember having understood them before, but
since we are content with this blind thought and do not pursue the
resolution of notions far enough, it happens that a contradiction
involved in a very complex notion is concealed from us'' (Leibniz as
translated by RA in \cite[note\;13]{Ra20}).  The readers can judge for
themselves how compelling this evidence is for RA's claim.}
but only non-contradictory concepts can refer to possibly existing
things; i.e., only non-contradictory concepts can be regarded as
defined semantically.

RA claim that Leibniz considers infinitely large or infinitely small
quantities as contradictory.  They write:
\begin{quote}
But, as we have shown, Leibniz always claimed that infinite entities,
be they infinitely large or infinitely small, could not be considered
as genuine quantities without violating a constitutive property of
quantities given by the part-whole axiom.  Hence, they cannot be
introduced into the system without contradiction.
\cite[pp.\;433--434]{Ra20}
\end{quote}
RA go on to conclude that Leibniz could not have understood
incomparables in non-Archi\-medean terms.

However, as we noted in connection with the \emph{infinita terminata},
RA provide no evidence that Leibniz considered all infinities
(including bounded ones) as contradictory.  Therefore, RA's conclusion
is based on an erroneous premise and a conflation of the notions of
bounded and unbounded infinity; see Section~\ref{s23} for details on
the distinction.  Section~\ref{s24} proposes an interpretation of the
distinction in modern mathematics.

While RA base their claims concerning the nature of Leibnizian
infinitesimals on their perception of alleged contradictions therein,
it is important to note that they furnish no corroborating evidence
whatsoever based on any textual analysis of Leibniz's actual
definitions of infinitesimals and infinite magnitudes; instead, they
rely on their claim of alleged contradictions.  In a passage we quoted
in our epigraph, RA even seem to admit that no such textual evidence
exists.

RA assert that Leibniz did not provide a semantic definition of
infinitesimals.  Could one expect him to provide such a definition?
Modern mathematical logic makes a clear distinction between syntactic
and semantic notions and between theory and model.  Modern philosophy
makes a distinction between procedures and ontology.  Such insights
are not easily attributable to Leibniz, though he did emphasize that
applying infinitesimals in geometry and physics should be independent
of metaphysical investigations.  Leibniz only specified the properties
he expected his inassignables to have, rather than providing models
that would allow him to ``have'' them in any semantic sense.  

As far as his ``theory of knowledge'' is concerned, Leibniz
specifically stated that one needn't get involved in metaphysical
questions as to the reality of such entities, and that the geometer's
task is limited to exploring the consequences of assuming the
existence of such mental constructs.  When he asserts that
infinitesimals ``have their proof with them'' what he seems to mean is
that their fruitful application
%
%
justifies their use, but not that they are metaphysically real (in a
few instances Leibniz specifically denies that material infinitesimals
exist, as in a letter to des Bosses).  The way Leibniz \emph{employs}
his inassignables is our only aid in interpreting LC.  The evidence
that he views them as consistent comes from the fact that he exploits
them in his mathematics.  Leibniz insisted many times on the
rigorousness of his calculus, and couldn't afford to exploit
inconsistent notions.

\subsection{RA \emph{vs} Knobloch}

RA repeatedly cite Knobloch's article (\cite{Kn02}, 2002) (without
providing either quotation or page number), as when they claim that
``the nub of the proof is an exploitation of the Archime\-dean
property to prove that quantities whose difference can be reduced to a
quantity smaller than any given quantity are equal''
\cite[p.\;414]{Ra20}.  In his earlier papers, Knobloch set out his
position clearly as follows (see e.g., \cite{Kn90}, 1990, p.\;42;
\cite{Kn94}, 1994, p.\;266--267).  Leibniz contrasts bounded infinity
and unbounded infinity.  Bounded infinity is fictional.  Unbounded
infinity exists in the physical world (and is even described by
Leibniz as ``actual''), but contradicts the part-whole axiom if taken
as a whole.  No claim concerning the presence of contradictions is
made with regard to bounded infinity, of which infinitesimals are one
of the manifestations.%
\footnote{Similarly, Jesseph notes that ``a central task for Leibniz's
fictionalism about infinitesimal magnitudes is to show that their
introduction is essentially harmless in the sense that it \emph{does
not yield a contradiction}'' (\cite[p.\;196]{Je15}; emphasis added).}
Knobloch's position in his early papers is in stark contrast with what
RA claim already in their abstract, namely that \emph{both} (i)
infinitesimals and (ii) unbounded infinity taken as a whole, lead to
contradictions.

\subsection
{RA's case against Bernoullian continua}

RA criticize the use of the term
\emph{Bernoullian continuum} by their opponents.  RA quote Leibniz as
follows: 
\begin{quote}
``We can conceive an infinite series consisting merely of finite
terms, or terms ordered in a decreasing geometric progression.  I
concede the infinite multiplicity of terms, but this multiplicity
forms neither a number nor one whole'' (To Johann Bernoulli, Feb. 24,
1699, GM III 575; A III 8, 66).  (Leibniz as translated by RA in
\cite[p.\;410]{Ra20})
\end{quote}
Here Leibniz clearly rejects infinite wholes.  RA then claim the
following: 
\begin{quote}
So we can see that the idea that Leibniz subscribed to a `Bernoullian
continuum', containing actual infinitesimals as its elements, is
emphatically rejected by Leibniz himself, writing to the same
Bernoulli.  \cite[p.\;411]{Ra20}
\end{quote}
However, RA provide no evidence that their opponents asserted that
Leibniz believed in infinite wholes.  Furthermore, the articles where
the Archimedean \emph{vs} Bernoullian distinction was introduced state
explicitly that they use Bernoulli's name merely because Bernoulli
exploited only B-track methods whereas Leibniz used both A-track and
B-track methods.  RA provide no evidence that infinitesimals violate
the part-whole axiom.

\subsection{What about Archimedean continua?}

The choice of target of RA's criticism is significant.  RA criticize
the use only of the term \emph{Bernoullian continuum} but not of the
term \emph{Archimedean continuum}.  Such an attitude may be related to
their undergraduate mathematical training resulting in a presentist
belief that `Archimedean' is the only natural way of thinking about
the continuum, and is symptomatic of butterfly-model commitments (see
Section~\ref{s17b}).  Meanwhile, Leibniz himself viewed any type of
infinite whole as contradicting the part-whole axiom.%
\footnote{\label{f26}Levey claims that ``For Leibniz there is no
genuine infinite \emph{magnitude}, great or small: no infinite line,
no infinite quantity, no infinite number'' (Levey \cite{Le21}, 2021,
p.\;147).  He seeks support for his claim in the following Leibnizian
passage from \cite{Le04}: ``It is perfectly correct to say that there
is an infinity of things, i.e.  that there are always more of them
than one can specify.  But it is easy to demonstrate that there is no
infinite number, nor any infinite line or other infinite quantity,
\emph{if} these are taken to be genuine wholes'' (emphasis added).
Levey appears not to have noticed that the quoted claim is
\emph{conditional} upon taking infinity to be a whole.  Unbounded
infinity taken as a whole contradicts the part-whole axiom; bounded
infinity doesn't (see Section~\ref{s23}).}

\subsection{Compactifications}
\label{s38}

RA misinterpret Leibnizian bounded infinities when they compare them
to compactifications \cite[p.\;421]{Ra20} in modern mathematics (such
as adding a point at infinity to an open infinite space).  Leibnizian
bounded infinities are an arithmetical extension rather than a
topological compactification.  Just as infinitesimals can be further
divided, according to Leibniz, \emph{infinita terminata} can be
expanded further.  For an infinite line from~$0$ to~$\mu$ there is
also a line from~$0$ to~$2\mu$, and so on.  While this would not make
any sense with compactifications, it fits well with the arithmetic
employing a nonstandard~$\mu$ (cf.\;Section~\ref{s24}).  Leibnizian
comparisons of infinitesimals to imaginaries and to ideal points of
intersection of parallel lines only refer to their status as fictions
possessing no instantiation in \emph{rerum natura}.

\subsection{RA on Kunen and paradoxes}

RA pursue a puzzling venture into modern set theory and the paradox of
``the set of all sets.''  They provide a lengthy quotation from
Kunen's textbook \cite{Ku80}, and propose the following comparison:
\begin{quote}
The standard set theorist can {\ldots} claim at the same time that
proper classes do not exist (in her axiomatic system), and that her
surface language remains neutral as regard this question of existence,
since there are other ways of interpreting it {\ldots} As we will see,
this is \emph{exactly the position endorsed by Leibniz} in his public
declarations regarding infinitesimals.  (\cite[note\;19]{Ra20};
emphasis added)
\end{quote}
Here RA appear to be comparing the fictionality of Leibnizian
infinitesimals with the fictionality of the set of all sets in modern
set theories.  However, such set-theoretic paradoxes involve entities
more similar to Leibnizian unbounded infinities than to his
infinitesimals%
\footnote{Arthur's conjured-up Leibniz specifically compares the
violation of the part-whole axiom by unbounded infinities, to the
set-theoretic paradox of the set of all ordinals, in
\cite[p.\;104]{Ar19b}.}
(see Section~\ref{s23}).  We therefore question the exactness of RA's
claim that ``this is exactly the position endorsed by Leibniz''
regarding infinitesimals.

See also Bl{\aa}sj\"o's analysis of the RA text in \cite{Bl20}.

\section
{Leibnizian calculus via modern infinitesimals}
\label{s24}

Recall that the canonical Zermelo--Fraenkel set theory (ZF) is a set
theory in the~$\in$-language.  Here~$\in$ is the two-place membership
relation.  Since the emergence of axiomatic nonstandard set theories
in the work of Hrbacek (\cite{Hr78}, 1978) and Nelson (\cite{Ne77},
1977), it became clear that it is possible to develop analysis in set
theories exploiting the more versatile~\st-$\in$-language.  Here {\st}
is a one-place predicate;~$\st(x)$ means ``$x$ is standard.''  A
recently developed theory SPOT (\cite{Hr20}, 2020) in the
~\st-$\in$-language is a conservative extension of ZF.  The theory
SPOT is a subtheory of theories developed by Hrbacek and Nelson.  For
a comprehensive treatment see Kanovei--Reeken (\cite{Ka04}, 2004) and
the references therein), as well as the survey by Fletcher et
al.\;(\cite{17f}, 2017).

Objections have been raised in the literature based on the assumption
that modern infinitesimal analysis allegedly depends on non-effective
foundational material such as the Axiom of Choice.  Thus, Ehrlich
comments:
\begin{quote}
The aforementioned ultrapower construction of hyperreal number systems
is another source of nonuniqueness since the construction depends on
an arbitrary choice of a nonprincipal ultrafilter.
\cite[note\;24, p.\;523]{Eh21}
\end{quote}
In a similar vein, Henle comments:
\begin{quote}
There was hope, when Abraham Robinson developed nonstandard analysis
{\ldots}, that intuition and rigor had at last joined hands.  His work
indeed gave infinitesimals a foundation as members of the set of
hyperreal numbers.  But it was an awkward foundation, dependent on the
Axiom of Choice.  {\ldots} Nonstandard analysis requires a substantial
investment (mathematical logic and the Axiom of Choice) but pays great
dividends.  \cite[pp.\;67, 72]{He99}
\end{quote}
There are many such comments in the literature.%
\footnote{Thus, V\"ath comments: ``Without the existence
of~$\delta$-free ultrafilters {\ldots} we were not able to construct
nonstandard embeddings.  {\ldots} [I]n the author's opinion {\ldots}
nonstandard analysis is not a good model for `real-world' phenomena''
\cite[p.\;85]{Va07}.  Easwaran and Towsner claim to ``point out
serious problems for the use of the hyperreals (and other entities
whose existence is proven only using the Axiom of Choice) in
describing the physical world in a real way'' \cite[p.\;1]{ET19}.
Sanders catalogs many such comments in \cite[p.\;459]{Sa20}.  A
rebuttal of the claims by Easwaran--Towsner and V\"ath appears in
\cite{19d}.  A rebuttal of related criticisms by Pruss \cite{Pr18b}
appears in \cite{21b}, \cite{21c}.  The theories developed in
\cite{Hr20} expose as factually incorrect a claim by Alain Connes to
the effect that ``as soon as you have a non-standard number, you get a
non-measurable set'' \cite[p.\;26]{Co07}.}
Such an assumption turns out to be incorrect.  Infinitesimal analysis
can be done in an axiomatic framework conservative over ZF, as shown
in \cite{Hr20}.  In particular, it can be developed without assuming
either the existence of nonprincipal ultrafilters or the axiom of
choice.

In this section, we will detail a formalisation of the procedures of
the Leibnizian calculus in SPOT.

\subsection{Assignable \emph{vs} inassignable} 

The predicate {\st} provides a formalisation of the Leibnizian
distinction between assignable and inassign\-able quantities, as used
in e.g., \cite{Le01c} (for an analysis see Bos \cite{Bo74} as well as
Section~\ref{s6} below).  Here \emph{assignable} corresponds to
\emph{standard} whereas \emph{inassignable} corresponds to
\emph{nonstandard}.  Leibniz distinguished notationally between the
inassignable~$dx$, $dy$ and the assignable~$(d)x, (d)y$; see
Section~\ref{s42}.

\subsection{Law of continuity and transfer}
\label{s32}

By the transfer principle, if~$\phi$ is an~$\in$-formula with standard
parameters, then the following entailment holds:
\[
\forall^{\st} x\; \phi(x) \to \forall x\; \phi(x).
\]
Here~$\forall^{\st}$ denotes quantification over standard entities
only.  SPOT's transfer principle provides a formalisation of the
Leibnizian law of continuity, to the effect that ``the rules of the
finite are found to succeed in the infinite'' (Robinson \cite{Ro66},
1966, p.\;266).  Namely, the rules of arithmetic for standard natural
numbers (finite numbers) extend to the rules of arithmetic of (all)
natural numbers; here the nonstandard ones correspond to the
Leibnizian \emph{infinita terminata} (see Section~\ref{s23}).  The
transfer principle can be thought of as a (belated) response to
Rolle's critique of the Leibnizian calculus and the law of continuity;
see \cite[Sections\;2.7--2.11]{18a} for details.%
\footnote{\label{f30}In connection with the law of continuity, Levey
observes: ``It is not hard to think of counterexamples to the law
[``the rules of the finite are found to succeed in the infinite''] and
at a minimum the law of continuity requires additional sharpening and
guidance to be used correctly in specific cases''
\cite[pp.\;153]{Le21}.  However, Levey fails to point out that such
sharpening and guidance are indeed provided in Robinson's classic
framework for infinitesimal analysis and its modern axiomatic versions
(see the beginning of the current Section~\ref{s24}).  Unlike Probst
\cite{Pr18}, Levey provides no hint of the existence of an alternative
to the syncategorematic reading.}

\subsection{Infinitesimals and Euclid V\!.4} 
\label{s43b}

A real number~$\epsilon$ is \emph{infinitesimal} if it satisfies the
formula
\[
\forall{n}\in\N\left[\st(n)\to|\epsilon|<\tfrac1n\right].
\]
Such a number necessarily satisfies~$\neg\st(\epsilon)$ if
$\epsilon\not=0$.  The formula above formalizes the violation of
Euclid Definition\;V.4 (see Section~\ref{V4}).

A real number smaller in absolute value than some standard real number
is called \emph{finite}, and otherwise \emph{infinite}.  Every
nonstandard natural number is infinite.

\subsection{Bounded \emph{vs} unbounded infinity: a formalisation}
\label{s43}

The theory SPOT enables a formalisation of the Leibnizian distinction
between bounded infinity and unbounded infinity (see
Section~\ref{s23}) as follows.  Unbounded infinity is exemplified by
the natural numbers~$\N$, whereas bounded infinity is exemplified
by~$\mu$ for a nonstandard (in Leibnizian terminology,
\emph{inassignable}) element~$\mu\in\N$ (thus we have
$\neg\st(\mu)$).%
\footnote{Bassler claims the following: ``while it is difficult if not
impossible to imagine a line starting at a particular point, going on
forever, and then terminating at another point, the idea that two
points on a line could be infinitely close to each other seems
considerably more palatable.  {\ldots} At any rate, this is much
different from Leibniz's conception that the model for the
mathematical infinite is the sequence of natural numbers, which has a
beginning but no end'' \cite[p.\;145]{Ba08}.  Arguably the two
conceptions are not ``much different'' since the interval
$[0,\mu]\subseteq\R$ would formalize Leibniz's bounded infinite line.
It is therefore not difficult to imagine a line starting at a
particular point, terminating at another point, and having infinite
length.}
In modern mathematics, the existence of~$\N$ as a whole depends on
infinitary axioms of set theory, which certainly contradict the
part-whole principle (viewed as axiomatic by Leibniz).  Meanwhile,
using~$\mu$ itself in the context of Skolem's model \cite{Sk33},
\cite{Sk34} of Peano Arithmetic entails no such foundational woes,
since Peano Arithmetic and the extended Zermelo-like theory of finite
sets are definitionally equivalent \cite[p.\;225]{Ta87}.%
\footnote{Fragments of nonstandard arithmetic are studied by Avigad
\cite{Av05}, Sommer and Suppes \cite{SS}, Nelson \cite{N}, Sanders
\cite{Sa20}, van den Berg and Sanders \cite{Va19}, Yokoyama
\cite{Yo10}, and others.}


\subsection{Discarding negligible terms}
\label{s35}

Leibniz often exploited the technique of discarding negligible terms
(see note~\ref{f21b}).  His \emph{transcendental law of homogeneity}
in \cite{Le10b} is a codification of the procedure; see \cite{12e} for
an analysis.  The technique is formalized in terms of the standard
part principle, to the effect that if~$|x|<r$ for some standard~$r$,
then there exists a standard~$x_0$ such that~$x-x_0$ is infinitesimal.

\subsection{Curves as infinilateral polygons}
\label{s46}

A curve~$\alpha(t)$,~$t\in[0,1]$ can be approximated by a polygon with
a nonstandard (inassignable) number~$\mu$ of sides, with vertices at
the points~$\alpha(\frac{i}{\mu})$ as~$i$ runs from~$0$ to~$\mu$.
Such approximations are sufficient for calculating the usual geometric
entities such as the tangent line and the osculating circle; see
Section~\ref{s28}.  For the corresponding treatment in SIA, see
Sections~\ref{s25} and \ref{s31}.

\subsection{Parabola as \emph{status transitus} of ellipses} 
\label{s37}

Leibniz notes that a parabola can be obtained as \emph{status
transitus} from a family of ellipses (see Section~\ref{s51}).
Consider a family of ellipses~$(E_t)$ such that every ellipse passes
through the origin and has one focus at~$(0,1)$ and the other at
$(0,t+1)$.  The family is parametrized by the distance~$t$ between the
foci.  Then the family includes an internal ellipse~$E_\mu$ for an
infinite value~$t=\mu$.  Taking standard part of the coefficients of a
normalized equation of $E_\mu$, we obtain the equation of a standard
parabola.

Other ways of formalizing Leibnizian principles are explored by Forti
(\cite{Fo18}, 2018).

\section{Infinitesimals of Smooth Infinitesimal Analysis}
\label{s3}

In this section we will analyze more closely the details of the
foundational approach adopted in {\law}, and their ramifications for
Arthur's comparison of Leibnizian calculus with {\law}.

\subsection
{Fictions and equivocations}
\label{s31}

The intuitionistic logic relied upon in {\law} enables an infinitesimal
$\epsilon$ to satisfy
\begin{equation}
\label{e16}
\neg\neg(\epsilon=0),
\end{equation}
or in words: ``$\epsilon$\;is not nonzero'' (see e.g., Bell
\cite{Be08}, 2008, formula\;(8.1), p.\;105).  On the other hand, one
can prove neither~$(\epsilon=0)$ nor\,~$\neg(\epsilon=0)$.

In this sense it could perhaps be said that the nilpotent
infinitesimals of {\law} hover at the boundary of existence.  Such a
mode of existence could be loosely described as fictional (or
``finessing of a commitment to the existence of infinitesimals'' as
Arthur puts it; see Section~\ref{s25}).  Even if one does use the term
\emph{fiction} to describe them, are such SIA fictions related to the
logical fictions (see Section~\ref{s13}) as posited by the
Ishiguro--Arthur (IA) syncategorematic reading of \emph{Leibnizian}
infinitesimals?

On the IA reading, an infinitesimal~$\epsilon$ is understood as
generated by ordinary real values, say~$\epsilon_n$, chosen small
enough to defeat a pre-assigned error bound.%
\footnote{Thus, speaking of Leibnizian infinite-sided polygons, Arthur
claims that ``this means that a curve can be construed as an ideal
limit of a sequence of such polygons, so that its length~$L$ will be
the limit of a sequence of sums~$ns$ of their sides~$s$ as their
number~$n\to\infty$'' (\cite[note 4, p.\;393]{Ar01}).}
But then the quantity~$H$ generated by~$H_n=\frac{1}{\epsilon_n}$ will
surely represent the inverse of the infinitesimal quantity~$\epsilon$,
so that we have a nonzero product \mbox{$\epsilon{H}=1$}.  Thus such
an~$\epsilon$ would necessarily be invertible and therefore provably
nonzero: \mbox{$\neg(\epsilon=0)$}, contradicting the {\law}
assumption \eqref{e16} that~$\neg\neg(\epsilon=0)$; cf.\;(Bell
\cite{Be08}, 2008, formula\;(8.2), p.\;105).  Similar reasoning
applies if Leibnizian~$\epsilon$ is taken to be a variable quantity
taking ordinary real values, inevitably leading to invertibility.

Thus Leibnizian infinitesimals, especially on the IA reading, are
unlike {\law} infinitesimals.  Arthur's claim of similarity between
them (see Section~\ref{s25}) is largely rhetorical, and amounts to
thinly veiled equivocation on the meaning of the qualifier
\emph{fictional}.

\subsection{The punctiform issue}

An important issue in the history of analysis concerns the nature of
the \emph{continuum} as being punctiform or nonpunctiform; Leibniz
endorsed the latter view.  In an interesting twist, Arthur (following
Bell \cite{Be08}, 2008, p.\;3) speaks of \emph{infinitesimals} as
being nonpunctiform, as noted by Knobloch; see Section\;\ref{s25}.
This tends to obscure the possibility that the exact nature of any
continuum Leibniz may have envisioned may be irrelevant to the actual
\emph{procedures} of his infinitesimal calculus.  Whether points
merely mark locations on the continuum (as was the case for Leibniz)
or are more fundamental to the actual make-up of the continuum (as in
modern set-theoretic approaches whether of Archimedean or
non-Archimedean type) is an issue mainly of foundational ontology that
is transverse to Leibniz's mathematical practice.%
\footnote{For an analysis of the procedures/ontology distinction see
B{\l}aszczyk et al.\;\cite{17d}.}

Leibniz may never have envisioned a punctiform continuum, but
infinitesimals that are \emph{not nonzero but not provably zero} (see
Section~\ref{s31}) have no known source in Leibniz's writings, either.%
\footnote{See Section~\ref{s21} for Arthur's translation of Leibniz's
endorsement of the law of excluded middle.}

\subsection{Variables and static quantities}
\label{s33}

One of Arthur's claimed similarities between LC and {\law} is ``the
dependence on variables (as opposed to the static quantities of both
Standard and Non-standard Analysis)'' (see Section~\ref{s17} for the
full quotation).  Let us examine Arthur's static claim.

The analog of the real line $\R$ of traditional analysis is the {\law}
line $R$.  Bell defines a part~$\Delta$ of the {\law} line~$R$ by the
formula~$\Delta=\{x\colon x^2=0\}$; see \cite[p.\;20]{Be08}. The
part~$\Delta$ does not reduce to~$\{0\}$ as discussed in our
Section~\ref{s31}.  The principle of microaffineness asserts the
following:
\begin{equation}
\label{e32}
\forall g\colon\Delta\to R \;\; \exists!\,b\in R \;\; \forall \epsilon
\in \Delta, \quad g(\epsilon) = g(0) + b.\epsilon
\end{equation}
(see \cite[p.\;21]{Be08}).  The fact that the equality is required
\emph{for\;all}~$\epsilon$ in~$\Delta$ is crucial for the uniqueness
of $b$.

Arthur notes that ``[t]he letter~$\epsilon$ then denotes a variable
ranging over~$\Delta$'' \cite[p.\;537]{Ar13}.  This is the only
mention of \emph{variables} in Arthur's discussion of the {\law} approach
to the derivative.  Concerning {\law} infinitesimals, Arthur writes:

\begin{quote}
The sense in which they are fictions in {\law}, however, is that
although it is denied that an infinitesimal neighbourhood of a given
point, such as~$0$, reduces to zero, it cannot be inferred from this
that there exists any point in the infinitesimal neighbourhood
distinct from~$0$.  (Arthur \cite{Ar13}, 2013, pp.\;572--573)
\end{quote}
Apparently the claimed similarity between LC and {\law} based on
`variables' amounts to the fact that the existence of a unique~$b$ as
in formula~\eqref{e32} depends crucially on universal
quantification~$\forall \epsilon$ over the part~$\Delta$.

What are we to make of Arthur's claimed contrast between such
variables and what he describes as ``the static quantities of both
Standard and Non-standard Analysis''?  It is true that in the
classical setting, if a function~$f$ is differentiable, the
derivative~$L=f'(x)$ can be computed from a \emph{single}
infinitesimal~$\epsilon\not=0$ by taking the standard part.
Namely,~$L$ is the standard part of
$\frac{f(x+\epsilon)-f(x)}{\epsilon}$; see Section~\ref{s28} and
Keisler (\cite{Ke86}, 1986).  Here universal quantification
over~$\epsilon$ is not needed.  Such an~$\epsilon$ could possibly be
described as ``static'' (or more precisely \emph{fixed}).

On the other hand, \emph{defining} differentiability in the classical
setting does require universal quantification over the
infinitesimal~$\epsilon$.  Namely, denoting by~$\forall^{\text{in}}$
universal quantification over nonzero infinitesimals and
by~$\exists_0^{\text{in}}$ existential quantification over (possibly
zero) infinitesimals, we have the following:
\[
f'(x)=L \text{\; if and only if \;}
\forall^{\text{in}}\epsilon\;\exists^{\text{in}}_0\ell
\left[\frac{f(x+\epsilon)-f(x)}{\epsilon}=L+\ell\right].
\]
In this sense, Robinsonian infinitesimals are no more ``static'' (in
an Arthurian sense) than the Lawvere--Kock--Bell ones.  Arthur's
rhetorical claim of similarity between Leibniz and {\law} on account
of ``variables'' and statics hinges on equivocation on the meaning of
the term \emph{variable}.

\section
{Law of continuity and \emph{status transitus} in \emph{Cum
Prodiisset}}
\label{s6}

Leibniz's unpublished text \emph{Cum Prodiisset} dates from around
1701.  The most adequate translation of Leibniz's \emph{law of
  continuity} as it appeared in \emph{Cum Prodiisset} was given by
Child:
\begin{quote}
In any supposed [continuous] transition, ending in any terminus, it is
permissible to institute a general reasoning, in which the final
terminus may also be included.%
\footnote{In the original: ``Proposito quocunque transitu continuo in
aliquem terminum desinente, liceat raciocinationem communem
instituere, qua ultimus terminus comprehendatur'' (Leibniz
\cite{Le01c}, 1701, p.\;40).}
(Leibniz as translated in Child \cite[p.\;147]{Ch})
\end{quote}
Child used the noncommittal term \emph{terminus}.  Meanwhile, both Bos
\cite[p.\;56]{Bo74} and Arthur \cite[p.\;562]{Ar13} use the term
\emph{limit} (or \emph{limiting case}) in their translations.  Such a
translation risks being presentist in that it is suggestive of the
modern notion of limit, confirming Cajori's \emph{grafting} diagnosis
(see Section~\ref{s11}).%
\footnote{Spalt uses the German term \emph{Grenze}
\cite[p.\;111]{Sp15} which is more successful because it differs from
the technical term \emph{Grenzwert} for limit.}

\subsection{\emph{Status transitus}}
\label{s51}

In the formulation cited above, Leibniz used the expression
\emph{status terminus}.  However, what Leibniz is really getting at
here is clear from his development that follows the formulation of the
law in \emph{Cum Prodiisset}.  Namely, the key issue is actually what
Leibniz refers to as the \emph{status transitus}:
\begin{quote}
\ldots{} a state of transition [\emph{status transitus} in the
original] may be imagined, or one of \emph{evanescence}, in which
indeed there has not yet arisen exact equality or rest or parallelism,
but in which it is passing into such a state, that the difference is
less than any assignable quantity; (Leibniz as translated by Child in
\cite[p.\;149]{Ch}; emphasis on `evanescence' added)
\end{quote}
Leibniz goes on to provide some examples:
\begin{quote}
{\ldots} also that in this state there will still remain some
difference, some velocity, some angle, but in each case one that is
infinitely small; and the distance of the point of intersection, or
the variable focus, from the fixed focus will be infinitely great, and
the parabola may be included under the heading of an ellipse, etc.%
\footnote{Leibniz goes on to write: ``Truly it is very likely that
Archimedes, and one who seems so have surpassed him, Conon, found out
their wonderfully elegant theorems by the help of such ideas; these
theorems they completed with \emph{reductio ad absurdum} proofs, by
which they at the same time provided rigorous demonstrations and also
concealed their methods'' (ibid.).  Leibniz's reference to the
concealment of the direct methods (similar to Leibniz's own) by
Archimedes indicates the existence of an alternative (infinitesimal)
method to be concealed, contrary to the IA thesis.  See related
comments on translation into Archimedean terms in note~\ref{f35}.}
(ibid.)
\end{quote}
The \emph{status transitus} is essentially the evanescent stage, as
when, for example, one evaluates an equation of an ellipse at an
infinite value of the parameter to produce an ellipse with infinite
distance between the foci.  The law of continuity asserts the
legitimacy of such a procedure.  The law as stated here is closely
related to Leibniz's formulation of the law of continuity in
\cite[p.\;93--94]{Le02} as quoted by Robinson: 
\begin{quote}
The rules of the finite are found to succeed in the infinite
(\cite{Ro66}, 1966, p.\;266)
\end{quote}
%
%
(see Section~\ref{s32} for a modern formalisation).  Namely,
substituting an infinite value of the parameter produces a
\emph{status transitus} given by a legitimate conic/ellipse.  Then the
finite part of the conic/ellipse is indistinguishable from a parabola
%
%
(see Section~\ref{s37}).

\subsection{Which things are taken to be equal?}
\label{s42}

How does Arthur handle the Leibnizian law of continuity?  Immediately
following his translation of Leibniz's law of continuity using the
presentist ``limiting case,'' Arthur quotes \emph{Cum Prodiisset} as
follows:
\begin{quote}
Hence it may be seen that in all our differential calculus there is no
need to call equal those things that have an infinitely small
difference, but those things are taken as equal that have no
difference at all, \ldots{} (Leibniz as translated by Arthur in
\cite[p.\;563]{Ar13})
\end{quote}
Here we deliberately interrupted the Leibnizian passage at a comma,
for reasons that will become clear presently.  At first glance, it may
appear that Leibniz is contradicting what he wrote in many texts
including his published response to Nieuwentijt (Leibniz \cite{Le95b},
1695) as well as (Leibniz \cite{Le10b}, 1710), to the effect that he
is working with a generalized notion of equality up to a negligible
term (see e.g., Section~\ref{s28}).
%
%
Note however that the unusually long sentence in \emph{Cum Prodiisset}
does not stop with the apparent endorsement of exact equality, but
rather continues to list a number of conditions and qualifications.
One of the qualifications is that exact equality results only after
one performs suitable algebraic simplifications and omissions (such as
omitting~$dx$ where appropriate):
\begin{quote}
\ldots provided that the calculation is supposed to have been rendered
general, applying equally to the case where the difference is
something and to where it is zero; and only when \emph{the calculation
has been purged as far as possible through legitimate omissions} and
ratios of non-vanishing quantities until at last application is made
to the ultimate case, is the difference assumed to be zero.%
\footnote{Arthur's translation is rather unsatisfactory.  Child
translated this passage as follows: ``Hence, it may be seen that there
is no need in the whole of our differential calculus to say that those
things are equal which have a difference that is infinitely small, but
that those things can be taken as equal that have not any difference
at all, provided that the calculation is supposed to be general,
including both the cases in which there is a difference and in which
the difference is zero; and provided that the difference is not
assumed to be zero until the calculation is purged as far as is
possible by legitimate omissions, and reduced to ratios of
non-evanescent quantities, and we finally come to the point where we
apply our result to the ultimate case'' (Leibniz as translated by
Child in \cite[pp.\;151--152]{Ch}).}
(Leibniz as translated by Arthur in \cite[pp.\;563--564]{Ar13};
emphasis added)
\end{quote}
Thus, while the simplification of the differential ratio in the case
of the parabola~$x^2=ay$ results in~$\frac{2x+dx}{a}$, subsequently
discarding the term~$dx$ does result in an exact
equality~$\frac{(d)y}{(d)x}=\frac{2x}{a}$.%
\footnote{See note \ref{f21b} on discarding terms.}
The Leibnizian~$(d)x$ and~$(d)y$ were rendered~$\underline d x$
and~$\underline dy$ in Bos \cite{Bo74}.  The notation~$(d)x$ and
$(d)y$ refers to assignable quantities whose ratio is the modern
derivative.%
\footnote{\label{f44}Levey's page-and-a-half discussion of this
Leibnizian derivation manages to avoid mentioning the crucial
Leibnizian distinction between~$dx$ and~$(d)x$.  He concludes that
``in fact~$dx$ stands for a variable finite quantity, and its behavior
reflects precisely the fact---ensured by the continuity of the curve
AY---that the difference between the abscissas can be taken as small
as one wishes, all the way to zero'' \cite[p.\;152]{Le21}.  However,
the conclusion applies only to the~$(d)x$, not the~$dx$ as per Levey.}
RA quote the following Leibnizian passage:
\begin{quote}
But if we want to retain~$dx$ and~$dy$ in the calculation in such a
way that they denote non-vanishing quantities even in the ultimate
case, let~$(d)x$ be assumed to be any assignable straight line
whatever; and let the straight line which is to~$(d)x$ as~$y$
or~$\Xlone{}_1\!Y$ is to~$\Xlone T$ be called~$(d)y${\ldots}
\end{quote}
What has Leibniz accomplished here?  Denoting the assignable value of
the differential quotient by~$L$, we note that Leibniz chooses an
assignable~$(d)x$ and then \emph{defines} a new quantity~$(d)y$ to be
the product~$L\,(d)x$.  RA claim that the method using~$(d)x$
and~$(d)y$ ``differs from the first one in the sense that it does not
rely on vanishing quantities'' \cite[p.\;439]{Ra20}.  However, the
method most decidedly does depend on vanishing quantities.  The
Leibnizian~$(d)x$ and~$(d)y$ may all be assignable but that is not
where the substantive part of infinitesimal calculus is.  The
nontrivial work goes into the evaluation of~$L$ using the appropriate
``legitimate omissions'', whereas introducing a new quantity by
setting~$(d)y=L\,(d)x$ is of lesser mathematical import.  RA's claim
to the contrary stems from paying exaggerated attention to the
rhetorical content of the surface language of the Leibnizian text, at
the expense of the mathematical content.  Their claim that ``the
introduction of the Law of Continuity as a postulate {\ldots} is fully
compatible with the `syncategorematic' view'' \cite[p.\;404]{Ra20}
remains unsubstantiated.

\subsection
{Are inassignable quantities dispensable?}

RA claim to provide a ``new reading'' \cite[p.\;404]{Ra20} of the
mature texts by Leibniz on the foundations of the calculus, in
particular the famous text \emph{Cum Prodiisset}.  They claim that
Leibniz presents his strategy based on the Law of Continuity as being
provably rigorous according to the accepted standards in keeping with
the Archimedean axiom and that a recourse to inassignable quantities
is therefore avoidable.  They base their claim on the Leibnizian
statement that one can switch to assignable quantities~$(d)x$,~$(d)y$
that keep the same ratio as~$dx$,~$dy$.

However, as noted in Section~\ref{s42}, the nontrivial part of the
calculus consists in determining the assignable value of the
differential ratio, involving omission of negligible terms.  Therefore
we cannot agree with RA's assertion that the Leibnizian strategy based
on the Law of Continuity is independent of inassignable quantities.

Both in his discussion of the example of the parabola and in the
passage quoted by Arthur, Leibniz speaks of \emph{omitting} terms.  To
elaborate, Leibniz specifically speaks of a stage where ``the
calculation has been purged as far as possible through legitimate
omissions.''  Namely, these are the stages in the calculation where
negligible terms are omitted.  Only afterwards does one attain a stage
where ``those things are taken as equal that have no difference at
all.''

This important feature of the calculations found in Leibniz does not
occur in a modern epsilon-delta definition of the concept of limit,
where~$\lim_{x\to c} f(x)=L$ is defined by the quantifier formula
\[
(\forall\epsilon>0)(\exists\delta>0)(\forall
x)\big[0<|x-c|<\delta\Longrightarrow |f(x)-f(c)|<\epsilon\big].
\]
No omission ever takes place here.  A freshman who attempts,
Leibniz-style, to omit terms on a test in a course following the
\emph{Epsilontik} approach will surely lose part of the credit.
Discarding terms is not only \emph{not} a feature of the modern
definition in terms of limits, but is on the contrary seen as a
typical freshman calculus error.  On the other hand, from the
viewpoint of Robinson's framework, Leibniz's procedure admits of
straightforward formalisation in terms of taking the standard part
(see Sections~\ref{s28} and~\ref{s35}).  Thus we cannot agree with
Arthur's appreciation:
\begin{quote}
It can be appreciated, I think, how close [Leibniz's approach] is to a
\emph{modern} justification of differentiation in terms of
\emph{limits}.  (Arthur \cite{Ar13}, 2013, p.\;564; emphasis added)
\end{quote}
Arthur's appreciation is yet another case of the grafting of the
modern theory of limits, in a mandatory Archimedean context
(symptomatic of butterfly-model thinking; see Section~\ref{s17b}), on
the calculus of Leibniz.

\section{Conclusion}

We have examined three interpretations of the procedures of the
Leibnizian calculus in the current literature, related respectively to
frameworks developed by Weierstrass, Lawvere, and Robinson.  Arthur
pursued comparisons of the Leibnizian procedures to the first two
frameworks.  We have argued that both of Arthur's readings of the
procedures of the Leibnizian calculus are less successful than the
interpretation in terms of the procedures of Robinson's framework for
analysis with infinitesimals.

Rabouin and Arthur claim to ``show that by 1676 Leibniz had already
developed an interpretation from which he never wavered, according to
which infinitesimals, like infinite wholes, cannot be regarded as
existing because their concepts entail contradictions.''  Their
position is at odds with that elaborated by Eberhard Knobloch in a
number of early texts, concerning the distinction between bounded
infinity and unbounded infinity in Leibniz.  While it is true, as RA
claim, that Leibniz's fictionalist view as elaborated in \emph{De
Quadratura Arithmetica} is consistent with later views (e.g., those
expressed in the february 1702 letter to Varignon), not only did RA
not show that Leibniz held infinitesimals to be contradictory but in
fact RA appear to admit (in the passage quoted in our first epigraph)
that Leibniz never expressed such an alleged conviction.

Leibniz held that infinitesimals are \emph{unlike} infinite wholes, in
that infinitesimals pertain to bounded infinity (\emph{infinita
terminata}) whereas infinite wholes pertain to unbounded infinity and
contradict the part-whole principle.

RA appear to base their conclusions on an analysis of the Leibnizian
texts that dwells excessively on the rhetorical content of the surface
language, and takes insufficient note of the mathematical content.
This is apparent in their claim that the Leibnizian method exploiting
the law of continuity was independent of the use of infinitesimals.
As we have shown, the claim does not hold up mathematically.  The same
is true of Arthur's comparison of Leibnizian infinitangular polygons
with curves in Smooth Infinitesimal Analysis.  To answer Probst's
question (see Section~\ref{one}), the reading based on Robinson's
framework (or its axiomatic conceptualisations such as SPOT
\cite{Hr20}) is arguably more successful on a number of issues,
including infinitangular polygons and fictionality of infinitesimals.

\section*{Acknowledgments} 

We are grateful to Viktor Bl{\aa}sj\"o, El\'{\i}as Fuentes Guill\'en,
Karel Hrbacek, Vladimir Kanovei, Eberhard Knobloch, Sam Sanders, David
Schaps, and David Sherry for helpful comments.  The influence of
Hilton Kramer (1928--2012) is obvious.

\end{document}